\documentclass{amsart}

\usepackage{amssymb}
\usepackage{amsmath}

\newtheorem{teo}{Theorem}[section]
\newtheorem{lemma}[teo]{Lemma}
\newtheorem{prop}[teo]{Proposition}
\newtheorem{corollario}[teo]{Corollary}
\theoremstyle{definition}
\newtheorem{definiz}[teo]{Definition}
\newtheorem{example}[teo]{Example}
\newcommand{\R}{\mbox{${\mathbb R}$}}
\newcommand{\Q}{\mbox{${\mathbb Q}$}}

\newcommand{\C}{\mbox{${\mathbb C}$}}
\newcommand{\Pal}{\mbox{${\mathbb P}$}}

\newcommand{\N}{\mbox{${\mathbb N}$}}

\newcommand{\id}{\sf id}

\theoremstyle{remark}  \newtheorem{nota}[teo]{Remark}
\numberwithin{equation}{section}

\begin{document}

\title[The dynamics of maps tangent to the identity]{The dynamics of
maps tangent to the identity\\ and with non-vanishing index}

\author{Laura Molino}

\address{Dipartimento di
Matematica, Universit\`a degli Studi di Parma, Viale G. P. Usberti
53/A, I-43100, Parma, Italy.} \email{laura.molino@unipr.it}


\begin{abstract}
Let $f$ be a germ of holomorphic self-map of $\C^2$ at the origin
$O$ tangent to the identity, and with $O$ as a non-dicritical
isolated fixed point. A parabolic curve for $f$ is a holomorphic
$f$-invariant curve, with $O$ on the boundary, attracted by $O$
under the action of $f$. It has been shown in \cite{Ab} that if the
characteristic direction $[v]\in\Pal(T_O\C^2)$ has residual index
not belonging to $\Q^+$, then there exist parabolic curves for $f$
tangent to $[v]$. In this paper we prove, with a different method,
that the conclusion still holds just assuming that the residual
index is not vanishing (at least when $f$ is regular along $[v]$).
\end{abstract}

\maketitle

\section{Introduction}
One of the most interesting areas in the study of local dynamics in
several complex variables is the dynamics near the origin $O$ of
maps tangent to the identity, that is of germs of holomorphic
self-maps $f\colon\C^n\to\C^n$ such that $f(O)=O$ and
$df_O=\hbox{\id}$. When $n=1$ the dynamics is described by the
celebrated Leau-Fatou flower theorem; when $n>1$ we are still far
from understanding the complete picture, even though very important
results have been obtained in recent years (see, e.g., \cite{E},
\cite{W}, \cite{Ha}, \cite{Ab}).

In this paper we want to investigate conditions ensuring the
existence of parabolic curves (the 2-variables analogue of the
petals in the Leau-Fatou flower theorem) for maps tangent to the
identity in dimension~2. Let us first recall some definitions and
useful results concerning maps tangent to the identity. Let $f$ be a
germ of holomorphic self-map of $\C^2$ fixing the origin and tangent
to the identity; we can write $f=(f_1,f_2)$, and let
$f_j=z_j+P_{j,\nu_j}+P_{j,\nu_j+1}+\cdots$ be the homogeneous
expansion of $f_j$ in series of homogeneous polynomials, where $\deg
P_{j,k}=k$ (or $P_{j,k}\equiv 0$) and $P_{j,\nu_j}\not\equiv 0$. The
{\em order} $\nu(f)$ of $f$ is defined by $\nu(f)=\min
\{\nu_1,\nu_2\}$. We say that the the origin is {\em dicritical} if
we have $z_2P_{1,\nu(f)}(z_1,z_2)\equiv z_1P_{2,\nu(f)}(z_1,z_2)$.

A {\em parabolic curve} for $f$ at the origin is an injective
holomorphic map $\varphi\colon\Delta\to\C^2$ satisfying the
following properties:

(i) $\Delta$ is a simply connected domain in~$\C$, with
$0\in\partial\Delta$;

(ii) $\varphi$ is continuous at the origin, and $\varphi(0)=O$;

(iii) $\varphi(\Delta)$ is invariant under $f$, and
$(f|_{\varphi(\Delta)})^n\to O$ as $n\to\infty$.

\noindent Furthermore, if $[\varphi(\zeta)]\to[v]\in\Pal^1$ as
$\zeta\to0$, where $[\cdot]$ denotes the canonical projection of
$\C^2\setminus\{O\}$ onto $\Pal^1$, we say that $\varphi$ is {\em
tangent to $[v]$} at the origin.

A {\em characteristic direction} for $f$ is a point
$[v]=[v_1:v_2]\in \Pal^1$ such that there is $\lambda \in\C$ so that
$P_{j,\nu(f)}(v_1,v_2)=\lambda v_j$ for $j=1,2$. If $\lambda\neq 0$,
we say that $[v]$ is {\em nondegenerate}; otherwise, it is {\em
degenerate}. Characteristic directions arise naturally if we want to
investigate the existence of parabolic curves tangent to some
direction $[v]$. In fact, Hakim observed that if there exist
parabolic curves tangent to a direction $[v]$ then this direction is
necessarily characteristic (\cite{Ha}). However, Hakim was able to
prove the converse for nondegenerate characteristic directions only:

\begin{teo}[\'Ecalle, \cite{E}; Hakim, \cite{Ha}, \cite{Ha1}]
\label{Hakim} Let $f$ be a (germ of) holomorphic self-map of $\C^2$
fixing the origin and tangent to the identity. Then for every
nondegenerate characteristic direction $[v]$ of $f$ there are
$\nu(f)-1$ parabolic curves tangent to $[v]$ at the origin.
\end{teo}

When $f$ has {\em no} nondegenerate characteristic directions, this
theorem gives no information about the dynamics of $f$. Furthermore,
there are examples of parabolic curves tangent to degenerate
characteristic directions.

\begin{example}
Let us consider the germ $f$ given by
\begin{displaymath}
\begin{cases}
f_1(z,w)=z+ zw+ w^2- z^3 +O(z^2w,zw^2,w^3,z^4),\\
f_2(z,w)=w[1+z+w+O(z^2,zw,w^2)].
\end{cases}
\end{displaymath}
We observe that $[1:0]$ is a degenerate characteristic direction.
The line $\{w=0\}$ is $f$-invariant, and inside it $f$ acts as the
function $z-z^3+O(z^4)$. The classical Leau-Fatou theory then shows
that there exist $2$ parabolic curves for $f$ tangent to $[1:0]$ at
the origin (\cite{F}, \cite{L}).
\end{example}

A further step toward the understanding of the dynamics in an
neighbourhood of an isolated fixed point has been done by Abate, who
gave a complete generalization of the Leau-Fatou flower theorem in
$\C^2$:

\begin{teo}[Abate, \cite{Ab}]
Let $f$ be a (germ of) holomorphic self-map of $\C^2$ tangent to the
identity and such that the origin is an isolated fixed point. Then
there exist (at least) $\nu(f)-1$ parabolic curves for $f$ at the
origin.
\end{teo}

The proof of this theorem is based on the possibility of modifying
the geometry of the ambient space via a finite number of blow-ups,
and of defining a {\em residual index} ${\sf Ind}(\tilde{f},S,p)\in
\C$, where $\tilde{f}$ is a holomorphic self-map of a complex
2-manifold $M$ which is the identity on a 1-dimensional submanifold
$S$, and $p\in S$. It turns out that this index is either not
defined anywhere on $S$, in which case we say that $\tilde{f}$ is
{\em degenerate} along $S$ (or {\em non-tangential} to $S$ in the
terminology of \cite{ABT}, where it is described a far-reaching
approach to indices for holomorphic self-maps), or it is everywhere
defined, and then we say that $\tilde{f}$ is {\em nondegenerate}
along $S$ (respectively, {\em tangential} to $S$).

In particular, Abate gave a generalization of Theorem \ref{Hakim} to
those characteristic directions whose residual index is not a
non-negative rational number.

\begin{teo}[Abate, \cite{Ab}]\label{Abate}
Let $f$ be a (germ of) holomorphic self-map of $\C^2$ tangent to the
identity and such that the origin is an isolated fixed point. Let
$[v]$ be a characteristic direction of $f$ such that ${\sf
Ind}(\tilde{f},\Pal^1,[v]) \notin \Q^+$ (here $\Pal^1$ is the
exceptional divisor of the blow-up of the origin, and $\tilde{f}$ is
the blow-up of $f$). Then there are (at least) $\nu(f)-1$ parabolic
curves for $f$ tangent to $[v]$ at the origin.
\end{teo}

The theory about the existence of parabolic curves tangent to a
direction $[v]$, for maps tangent to the identity, is thus almost
complete, but there are still examples where the previous results
cannot be applied.

\begin{example}
Let us consider the map $f$ given by
\begin{displaymath}
\begin{cases}
f_1(z,w)=z+ zw+ O(w^2,z^3,z^2w),\\
f_2(z,w)=w+ 2w^2 + bz^3 + z^4+O(z^5,z^2w,zw^2,w^3),
\end{cases}
\end{displaymath}
with $b\neq0$. We observe that $[v]=[1:0]$ is a degenerate
characteristic direction for $f$ with ${\sf
Ind}(\tilde{f},\Pal^1,[v])=1$. Hence, we cannot say anything about
the dynamics of $f$ in the direction $[1:0]$ using Theorem
\ref{Hakim} or Theorem \ref{Abate}.
\end{example}

As a corollary of our work (see Corollary~1.7) we shall be able to
prove the existence of parabolic curves tangent to $[1:0]$ also for
this example (and many others). This will be a consequence of a more
general result; to state it precisely we need to recall another set
of definitions.

Any germ $g$ of holomorphic function defined in a neighbourhood of
the origin $O$ has a homogeneous expansion as an infinite sum of
homogeneous polynomials, $g=P_0+P_1+\cdots$, with $\deg P_j=j$ (or
$P_j\equiv 0$); the least $j\geq 0$ such that $P_j$ is not
identically zero is the {\em order} $\nu(g)$ of $g$. Let $f$ be a
germ of holomorphic self-map of $\C^2$ fixing the origin and tangent
to the identity. We can write $f_1=z+g$ and $f_2=w+h$. Set
$l=\gcd(g,h)$, and write $g=lg^o$ and $h=lh^o$. The {\em pure order}
of $f$ at the origin is $\nu_o(f,O)=\min \{ \nu(g^o),\nu(h^o)\}$. We
say that the origin is a {\em singular point}  for $f$ if
$\nu_o(f,O)\ge1$. If $\hbox{Fix}(f)$ has (at least) two (local)
components intersecting at the origin, we say that the origin is a
{\em corner}.

More generally, let $S$ be a subset of a complex 2-manifold $M$; we
shall denote by ${\sf End}(M,S)$ the set of germs about $S$ of
holomorphic self-maps of $M$ sending $S$ into itself. Let $f\in {\sf
End}(M,S)$ be such that $f|_S=\hbox{\id}_S$, where $S$ is a
1-dimensional submanifold of $M$, and assume that $df$ acts as the
identity on the normal bundle of $S$ in $M$. Then we can extend the
definitions of pure order, corner, singular point, and dicritical
point to any $p\in S$, simply by choosing a chart of $M$ centered at
$p$ and considering the local expression of $f$.

Our main theorem can now be stated as follows:

\begin{teo}\label{teomain}
Let $S$ be a 1-dimensional submanifold of a complex 2-manifold $M$
and let $f\in {\sf End}(M,S)$ be such that $f|_S=\hbox{\id}_S$.
Assume that $df$ acts as the identity on the normal bundle of $S$ in
$M$ and let $f$ be tangential to $S$. If  $p\in S$ is a singular
point of $f$, not a corner, with $\nu_o(f,p)=1$ and ${\sf
Ind}(f,S,p) \neq 0$ then there exist parabolic curves for $f$ in
$p$.
\end{teo}

An important application of this result is the following. Starting
from a map $f\in {\sf End}(\C^2,O)$ tangent to the identity and
blowing up the origin, we obtain an $\tilde{f}\in {\sf End}(M,S)$,
where $S\cong\Pal^1$ is the exceptional divisor of the blow-up. It
turns out that $\tilde f|_S=\hbox{\id}_S$, and that $d\tilde f$ acts
as the identity on the normal bundle of $S$ in $M$. Furthermore, if
the origin is an isolated fixed point of $f$, then no point $p\in S$
is a corner, and $\tilde f$ is tangential to $S$ if and only if the
origin is non-dicritical for $f$. If the origin is dicritical for
$f$, then all directions are characteristic, and there are parabolic
curves tangent to all but a finite number of them; so we concentrate
on the non-dicritical case. Every nondegenerate characteristic
direction of $f$ is a singular point for $\tilde f$, and every
singular point of $\tilde f$ is a characteristic direction of $f$.
Moreover, if $p\in S$ is not singular then no infinite orbit can get
arbitrarily close to $p$ (in particular, no infinite orbit can
converge to $p$, and thus there can be no parabolic curves at $p$);
therefore, from a dynamical point of view only singular points are
interesting (for the proof of all these assertions see \cite{Ab}).

Finally, we shall say that $f$ is {\em regular} along the
characteristic direction $[v]\in\Pal^1$ if the pure order of $\tilde
f$ at $[v]$ is one. This is just a technical condition almost always
satisfied (for instance, it is satisfied by the map in Example~1.5).
Here is the promised

\begin{corollario}\label{cormain}
Let $f\in {\sf End}(\C^2,O)$ be tangent to the identity with the
origin as a non-dicritical isolated fixed point. Let $[v]\in \Pal^1$
be a characteristic direction of $f$ and assume $f$ is regular along
$[v]$ with ${\sf Ind}(\tilde{f},\Pal^1,[v])\neq 0$ (here we identify
$\Pal^1$ with the exceptional divisor of the blow-up of the origin,
and $\tilde{f}$ is the blow-up of $f$). Then there exist parabolic
curves for $f$ tangent to $[v]$ at the origin.
\end{corollario}

\vskip0.2cm

 \noindent\textbf{Acknowledgements.}
I would like to thank M. Abate for his encouragement and support
during the preparation of this article. I also thank the referee for
useful comments and suggestions.

\section{Blow-up}

Since blow-ups will be a fundamental tool in our study we recall
some basic definitions, referring to \cite{Ab1} for details. Let $M$
be a complex 2-manifold, and let $p\in M$. The {\em blow-up} of $M$
in $p$ is the set $\tilde{M}=(M\setminus \{p\})\cup \Pal(T_pM)$,
together with the projection $\pi\colon\tilde{M} \to M$ given by
$\pi|_{M\setminus \{p\}}={\id}_{M\setminus \{p\}}$ and
$\pi|_{\Pal(T_pM)}\equiv p$. The set $S=\Pal(T_pM)=\pi^{-1}(p)$ is
the {\em exceptional divisor} of the blow-up.

It is possible to endow $\tilde{M}$ with a structure of
2-dimensional complex manifold. Fix a chart $\psi=(z_1,z_2)\colon
U\to \C^2$ of $M$ centered at $p$. For $j=1,2$, set
\begin{displaymath}
U_j=(U\setminus\{z_j=0\})\cup(S\setminus \hbox{Ker}(dz_j|_p)),
\end{displaymath}
and let $\chi_j: U_j\to \C^2$ be given by
\begin{displaymath}
\chi_j(q)_h=\left\{ \begin{array}{ll} z_j(q) & \textrm{if $j=h$ and
$q\in U\setminus\{z_j=0\}$},\\
\frac{z_h(q)}{z_j(q)} & \textrm{if $j\neq h$ and $q\in
U\setminus\{z_j=0\}$},\\
\frac{d(z_h)_p(q)}{d(z_j)_p(q)} & \textrm{if $j\neq h$ and
$q\in S\setminus \hbox{Ker}(dz_j|_p))$},\\
0 & \textrm{if $j=h$ and $q\in S\setminus \hbox{Ker}(dz_j|_p))$}.
\end{array} \right.
\end{displaymath}
Then the charts $(U_j,\chi_j)$, together with an atlas of
$M\setminus \{p\}$, endow $\tilde{M}$ with a structure of
2-dimensional complex manifold such that the projection $\pi$ is
holomorphic everywhere and given by
\begin{displaymath}
[\psi \circ \pi \circ \chi_j^{-1}(w)]_h = \left\{ \begin{array}{ll}
w_j & \textrm{if $j=h$},\\
w_jw_h & \textrm{if $j\neq h$}.
\end{array} \right.
\end{displaymath}
Let $f\in {\sf End}(M,p)$ be such that $df_p$ is invertible. Then
(see \cite{Ab1}) there exists a unique map $\tilde{f}\in {\sf
End}(\tilde{M},S)$, the {\em blow-up} of $f$ at $p$, such that $\pi
\circ \tilde{f}=f\circ \pi$. The action of $\tilde{f}$ on $S$ is
induced by the action of $df_p$ on $\Pal(T_pM)$; in particular, if
$df_p=\hbox{\id}$, then $\tilde{f}|_S=\hbox{\id}_S$. Finally, if $R$
is a submanifold of $M$, then the {\em proper transform} of $R$ is
$\tilde{R}=\overline{\pi^{-1}(R\setminus \{p\})}$.

\section{residual index}

Let us fix notations and definitions that we use in the paper.
$\mathcal{O}_2$ denotes the ring of germs of holomorphic functions
defined in a neighbourhood of the origin $O$ of $\C^2$. If $S$ is a
1-dimensional submanifold of a complex 2-manifold $M$, a chart
$(U,\varphi)$ of $M$ about $p\in S$ is {\em adapted} to $S$ if
$U\cap S=\varphi^{-1}(\{(z,w) | z=0\})$. Let $f\in {\sf End}(M,S)$
be such that $f|_{S}=\hbox{\id}_{S}$ and assume that $df$ acts as
the identity on the normal bundle of $S$ in $M$. Then in an adapted
chart $(U,\varphi)$ centered at $p\in S$ we can write
\begin{displaymath}
\begin{cases}
f_1(z,w)=z+ z^{\mu +2}\hat{g}(z,w),\\
f_2(z,w)=w+z^{\nu +1}\hat{h}(z,w),
\end{cases}
\end{displaymath}
for suitable $\hat{g}$, $\hat{h}\in \mathcal{O}_2$ and $\mu$, $\nu
\in \N\cup \{\infty\}$, where $\mu=\infty$ (resp., $\nu=\infty$)
means $\hat{g}\equiv 0$ (resp., $\hat{h}\equiv 0$), and where $z$
does not divide either $\hat{g}$ or $\hat{h}$. After having
introduced the (locally defined) meromorphic function
\begin{displaymath}
k(w)=\lim_{z\to 0}\frac{f_1(z,w)-z}{z(f_2(z,w)-w)}=\left\{
\begin{array}{ll} 0 & \textrm{if $\mu>\nu$},\\
(\hat{g}/\hat{h})(0,w) & \textrm{if $\mu=\nu$},\\
\infty & \textrm{if $\mu<\nu$},
\end{array} \right.
\end{displaymath}
Abate in \cite{Ab} called $f$ {\em degenerate} along $S$ if $k\equiv
\infty$. This definition is well posed (i.e., does not depend on the
adapted chart chosen); furthermore  he defined the {\em residual
index} of $f$ at $p$ {\em along} $S$ by setting
\begin{displaymath}
{\sf Ind}(f,S,p)=\hbox{Res}_0(k(w));
\end{displaymath}
again, it is independent of the adapted chart chosen, and we observe
that it might be nonzero at singular points of $f$ only. See
\cite{Ab} and \cite{ABT} (in particular the latter) for a more
thorough and deeper discussion of residual indices and related
topics.

\section{Dynamics}

In this section we give the proof of Theorem \ref{teomain} and
Corollary \ref{cormain}. First, we introduce some definitions. Let
us consider $f\in {\sf End}(\C^2,O)$ tangent to the identity, and
$[v]\in\Pal^1$ a characteristic direction for $f$.

\begin{definiz}
We say that $f$ is {\em regular} along $[v]$ if
$\nu_o(\tilde{f},[v])=1$, where $\tilde{f}$ is the blow-up of $f$ at
the origin.
\end{definiz}

\begin{nota}
If $[v]$ is a nondegenerate characteristic direction for $f$ then
$f$ is regular along $[v]$. Generally, the converse is not true.
\end{nota}

\begin{example}
Let us consider a map $f\in {\sf End}(\C^2,O)$ of order $2$,
\begin{displaymath}
\begin{cases}
f_1(z,w)=z+a_{2,0}z^2+a_{1,1}zw+a_{0,2}w^2+\cdots ,\\
f_2(z,w)=w+b_{2,0}z^2+b_{1,1}zw+b_{0,2}w^2+b_{3,0}z^3\cdots .
\end{cases}
\end{displaymath}
We observe that $[v]=[1:0]$ is a characteristic direction for $f$ if
and only if $b_{2,0}=0$. In this case, $f$ is regular along $[1:0]$
if it is nondegenerate or it is degenerate but $b_{3,0}\neq 0$.
\end{example}

We can finally start working. Let $S$ be a $1$-dimensional
submanifold of a complex $2$-manifold $M$, and let $f$ be a germ
about $S$ of holomorphic self-map of $M$ such that
$f|_S=\hbox{\id}_S$. Assume that $df$ acts as the identity on the
normal bundle of $S$ in $M$, and assume that $f$ is tangential to
$S$. Let $p\in S$ be a non corner and a singular point, with
$\nu_o(f,p)=\ 1$ and ${\sf Ind}(f,S,p)\neq 0$. We want to prove that
there exist parabolic curves for $f$ in $p$. We work in a chart
adapted to $S$ and centered at $p$. Then we can write
\begin{equation}\label{form}
\begin{cases}
f_1(z,w)=z + z^rA_1(z,w),\\
f_2(z,w)=w + z^rB_1(z,w),
\end{cases}
\end{equation}
for suitable $A_1$, $B_1\in \mathcal{O}_2$,  $r\in\N^*$, with
$\min\{\nu(A_1),\nu(B_1)\}=1$ and $\gcd(A_1,B_1)=1$. Since $f$ is
tangential to $S$, we must have $A_1(z,w)=zA_0(z,w)$ with $\nu
(A_0)\geq 0$. Let
\begin{displaymath}
\begin{cases}
A_0(z,w)=a_{0,0} + (a_{1,0}z + a_{0,1}w) + (a_{2,0}z^2 + a_{1,1}zw
+ a_{0,2}w^2) +\cdots,\\
B_1(z,w)=(b_{1,0}z + b_{0,1}w) + (b_{2,0}z^2 + b_{1,1}zw +
b_{0,2}w^2) + \cdots,
\end{cases}
\end{displaymath}
be the homogeneous expansions of $A_0$ and $B_1$ in series of
homogeneous polynomials. Since $\gcd(A_1,B_1)=1$, $z$ does not
divide $B_1(z,w)$; thus we have $b_{0,j}\neq 0$, for some $j\geq 1$.
An easy calculation shows that
\begin{equation}\label{index}
{\sf Ind}(f,S,p) = \textrm{Res} _0\left( \frac{A_0(0,w)}{B_1(0,w)}
\right)=\textrm{Res} _0\left(
\frac{a_{0,0}+a_{0,1}w+a_{0,2}w^2+\cdots}
{b_{0,1}w+b_{0,2}w^2+b_{0,3}w^3+\cdots}\right).
\end{equation}
Set
\begin{displaymath}
m:=\min\{h \in \N : a_{0,h}\neq 0\},
\end{displaymath}
\begin{displaymath}
n:=\min\{j \in \N : b_{0,j}\neq 0\}.
\end{displaymath}
We observe that it must be $m<n$ because otherwise, by
~\eqref{index}, ${\sf Ind}(f,S,p)=0$.

Let us dispose of the easier cases first.

\begin{teo}
Assume that either

(a) $m<n-1$, or

(b) $m=n-1$ and ${\sf Ind}(f,S,p)\neq n$, or

(c) $m=0$, $n=1$ and ${\sf Ind}(f,S,p) = 1$.

Then there exist (at least) $r+m(r+1)$ parabolic curves for $f$ at
the origin.
\end{teo}

\begin{proof}
Let $n=1$. By \eqref{index} it follows that ${\sf
Ind}(f,S,p)=a_{0,0}(b_{0,1})^{-1}$.

Assume that ${\sf Ind}(f,S,p)\neq 0,1$. Then $[v]:=\left[1:
b_{1,0}(a_{0,0}-b_{0,1})^{-1}\right]$ is a nondegenerate
characteristic direction of $f$. In fact, by \eqref{form}, $[1:w]$
is a characteristic direction if and only if
\begin{equation}\label{system}
\begin{cases}
a_{0,0}=\lambda ,\\
b_{1,0} + b_{0,1}w=\lambda w,
\end{cases}
\end{equation}
which is satisfied by $w=b_{1,0}(a_{0,0}-b_{0,1})^{-1}$ and with
$a_{0,0} \neq 0$. Hence, we have $r$ parabolic curves for $f$
tangent to $[v]$ at the origin (cf. \cite{Ha}, \cite{Ha1}).

Assume now ${\sf Ind}(f,S,p)= 1$. We can write \eqref{form} as
\begin{equation}
\nonumber
\begin{cases}
f_1(z,w)=z + a_{0,0}z^{r+1} + O(z^{r+2},z^{r+1}w),\\
f_2(z,w)=w[1 + a_{0,0}z^r + O(z^{r+1},z^rw)] + O(z^{r+1}),
\end{cases}
\end{equation}
with $a_{0,0}\neq 0$. Setting $Z=\alpha z$, where
$\alpha^r=-a_{0,0}$, we reduce to
\begin{displaymath}
\begin{cases}
f_1(z,w)=z - z^{r+1} + O(z^{r+2},z^{r+1}w),\\
f_2(z,w)=w[1 - z^r + O(z^{r+1},z^rw)] + O(z^{r+1}).
\end{cases}
\end{displaymath}
The existence of $r$ parabolic curves at the origin for a map of
this form is then a consequence of the results of \cite{Ha}, adapted
as in \cite{Ha1} if $r>1$. We observe that these curves are tangent
to $[0:1]$ if $b_{1,0}\neq 0$, otherwise to $[1:0]$.

Let $n>1$, and assume $m<n-1$, or otherwise $m=n-1$ and ${\sf
Ind}(f,S,p)\neq n$.

If $m=0$, $f$ has a nondegenerate characteristic direction. In fact,
the system defining characteristic directions becomes
\begin{displaymath}
\begin{cases}
a_{0,0}=\lambda,\\
b_{1,0}=\lambda w,
\end{cases}
\end{displaymath}
with $\lambda \neq 0$. Then, again from \cite{Ha} and \cite{Ha1}, we
have the existence of $r$ parabolic curves for $f$ tangent to
$[1:b_{1,0}/a_{0,0}]$ at the origin.

If $m\geq 1$ we know that
\begin{displaymath}
\left\{ \begin{array}{ll} a_{0,0}=\cdots =a_{0,m-1}=0, \qquad \qquad
\qquad a_{0,m}\neq 0,\\
b_{0,1}=\cdots =b_{0,m}=\cdots =b_{0,n-1}=0, \quad b_{0,n}\neq 0.
\end{array} \right.
\end{displaymath}
Let us consider the linear chain starting in $p\in S$ as defined in
\cite{CS}. It is a finite sequence of projective lines intersecting
each other transversally and at most in one point, obtained by a
finite sequence of blow-ups reducing corners to simpler singular
points. Finiteness, as the possibility to obtain corners in a
simpler form, is assured by the reduction theorem \cite{Ab}. For our
purposes it is sufficient to consider only a part of this linear
chain. Blowing-up $f$ in $p$, let us denote by $\tilde {f}^{[1]}$
its blow-up, by $S_1$ the exceptional divisor, and by $\tau(p)$ the
intersection point between $S_1$ and $\tilde{S}$, the proper
transform of $S$.

{\em Notations}. To avoid heavy notation, at any step of this
construction, we will continue to denote by $S$ the proper transform
$\tilde{S}$.

The second step consists in blowing-up $\tilde {f}^{[1]}$ in
$\tau(p)$ and denoting by $\tilde {f}^{[2]}$ its blow-up, by $S_2$
the exceptional divisor, and by $\tau(\tau(p))=\tau^{2}(p)$ the
intersection point between $S_2$ and $S$. For $k\leq m$, iterating
this process, we can define $\tilde {f}^{[k]}$ as the blow-up of
$\tilde {f}^{[k-1]}$ in $\tau^{k-1}(p)$, $S_k$ the exceptional
divisor and $\tau^{k}(p)$ the intersection point between $S_k$ and
$S$. We observe that this point is strictly related to the order of
the chain, as defined in \cite{CS}. By induction, it is easy to
check that, in a chart centered at $\tau^{k}(p)$, we have
\begin{displaymath}
\begin{cases}
\begin{aligned}
\tilde{f}^{[k]} _1(z,w)=&z + z^rw^{k(r+1)-1}[-kb_{1,0}z^2 +
(a_{0,m}-kb_{0,m+1})zw^{m-k+1}\\
&\qquad \qquad \qquad \quad \,+O(zw^2,zw^{m-k+2})],\\
\tilde{f}^{[k]} _2(z,w)=&w + z^rw^{k(r+1)-1}[b_{1,0}zw +
b_{0,m+1}w^{m-k+2} + O(zw^2,w^{m-k+3})].
\end{aligned}
\end{cases}
\end{displaymath}
Let $z_0:=[a_{0,m}-(m+1)b_{0,m+1}](m+1)^{-1}b_{1,0}^{-1}$. We remark
that $z_0\neq 0$. We claim that $[z_0:1]$ is a nondegenerate
characteristic direction for $\tilde{f}^{[m]}$. In fact, there
exists $\lambda\in\C\setminus \{0\}$ such that
\begin{equation}\label{system1}
\begin{cases}
-mb_{1,0}z_0^{r+2} + (a_{0,m}-mb_{0,m+1})z_0^{r+1}=\lambda z_0,\\
b_{1,0}z_0^{r+1}+ b_{0,m+1}z_0^r=\lambda;
\end{cases}
\end{equation}
an easy computation shows that it suffices to choose
\begin{displaymath}
\lambda =
\begin{cases}
b_{1,0}z_0^{r+1} \quad \textrm{if $m<n-1$},\\
\frac{a_{0,n-1}}{n}z_0^r \quad \textrm{if $m=n-1$}.
\end{cases}
\end{displaymath}
So, again from \cite{Ha} and \cite{Ha1}, there exist
$\nu(\tilde{f}^{[m]})-1=r+m(r+1)$ parabolic curves for
$\tilde{f}^{[m]}$ tangent to $[z_0:1]$ at the origin. Notice that a
parabolic curve for $\tilde{f}^{[m]}$ cannot intersect the
exceptional divisor since all points of the curve are attracted to
the origin. Therefore, the push-forward of a parabolic curve for
$\tilde{f}^{[m]}$ is a parabolic curve for $\tilde{f}^{[m-1]}$
tangent to $[0:1]$ at the origin. The iteration of this process will
give us $r+m(r+1)$ parabolic curves for $f$ tangent to $[0:1]$ at
the origin.
\end{proof}

\begin{nota}
If $m=n-1$ and ${\sf Ind}(f,S,p)=n$ then $z_0=0$, and $[z_0:1]$ is a
degenerate characteristic direction for $\tilde{f}^{[n-1]}$. We can
say much more: in fact $[z_0:1]$ coincides with $\tau^n(p)$, and it
is an irreducible singular point of type $(\star _2)$ for
$\tilde{f}^{[n]}$, and thus the order of the linear chain, starting
at $p$, is exactly $n$ (see \cite{Ab} and \cite{CS} for the
terminology). Therefore we cannot expect anything from the previous
construction, because it ends at such a point.
\end{nota}

Let us now deal with the hardest case. We shall prove the following:
\begin{teo}\label{m=n-1}
Let $n\geq 2$. If $m=n-1$ and ${\sf Ind}(f,S,p)=n$, then there
exists (at least) one parabolic curve for $f$ at the origin.
\end{teo}

First, we need another set of results.

Assume then $n\geq 2$, $m=n-1$ and ${\sf
Ind}(f,S,p)=\frac{a_{0,n-1}}{b_{0,n}}=n$. We can write \eqref{form}
as follows:
\begin{equation}\label{fstartn}
\begin{cases}
\begin{aligned}
f_1(z,w)&=z +a_{0,n-1} z^{r+1}w^{n-1} + a_{1,0}z^{r+2} +
O(z^{r+3},z^{r+2}w,z^{r+1}w^{n}),\\
f_2(z,w)&=w[1 + b_{0,n}z^rw^{n-1} + b_{1,1}z^{r+1} +
O(z^{r+2},z^{r+1}w,z^rw^n)]\\
& \quad \,\,\,\,+ b_{1,0}z^{r+1} + O(z^{r+2}),
\end{aligned}
\end{cases}
\end{equation}
with $r\in\N^*$ and $b_{1,0}\neq 0$, because $\nu_o(f,p)=1$.

In what follows we need to consider analytic changes of variables,
defined in a neighbourhood of the origin $O\in\mathbb{C}^2$, and
which involve roots and logarithms of the first complex variable
$z$. This makes sense if $z$ belongs to a simply connected open set
with $z=0$ on the boundary, for example, the $z$-plane with a cut.
In this case, we fix, once and for all, one branch of the roots and
the logarithm of $z$ in a way that $t=(\log z)^{\frac{1}{n}}$
belongs to an open sector of the $t$-plane with a vertex at the
origin and a positive central angle non exceeding $\frac{\pi}{n}$.

\begin{prop}
Let $(f_1,f_2)$ be a map of the form \eqref{fstartn}. Then we can
choose local coordinates, defined in a suitable simply connected
open set with $O$ on the boundary, relative to which the map takes
the form:
\begin{enumerate}
\item if $n=2$:
\begin{equation}\label{fhat}
\begin{cases}
\begin{aligned}
\hat{f}_1(z,w)&=z -a^{-1} z^{r+1}w -z^{r+1+\frac{1}{2}}(\log
z)^{\frac{1}{2}}\\
&\qquad +O\Big(z^{r+2}\log z,z^{r+1+\frac{1}{2}}w(\log
z)^{\frac{1}{2}},
z^{r+1}w^{2}\Big),\\
\hat{f}_2(z,w)&=w\bigg[1 -\frac{a^{-1}}{2}z^rw
-\frac{1}{2}z^{r+\frac{1}{2}}(\log
z)^{\frac{1}{2}} +\frac{1}{2} \frac{z^{r+\frac{1}{2}}}{(\log z)^{\frac{1}{2}}} \\
&\qquad + O\Big(z^{r+1}\log z,z^{r+\frac{1}{2}}w(\log
z)^{\frac{1}{2}},z^rw^2\Big)\bigg]  + z^{r+1+\frac{1}{2}}\psi_1(z),
\end{aligned}
\end{cases}
\end{equation}
for some $a\neq 0$;

\item if $n\geq 3$:
\begin{equation}\label{fhatn}
\begin{cases}
\begin{aligned}
\hat {f}_1(z,w)&=z-(n-1)a^{-1}z^{r+1+\frac{n-2}{n}}w(\log
z)^{\frac{n-2}{n}} - z^{r+1+\frac{n-1}{n}}(\log
z)^{\frac{n-1}{n}}\\
&\qquad +O\Big(z^{r+2}\log z,z^{r+1+\frac{n-1}{n}}w(\log
z)^{\frac{n-1}{n}} \Big)\\
&\qquad +O\Big(z^{r+1+\frac{n-3}{n}}w^{2}(\log
z)^{\frac{n-3}{n}},\ldots,z^{r+1}w^{n-1}\Big),\\
\hat{f}_2(z,w)&=w\bigg[1-\frac{n-1}{n}a^{-1}z^{r+\frac{n-2}{n}}w(\log
z)^{\frac{n-2}{n}} \\
&\qquad -\frac{1}{n}z^{r+\frac{n-1}{n}}(\log z)^{\frac{n-1}{n}}
+\frac{n-1}{n}\frac{z^{r+\frac{n-1}{n}}}{(\log z)^{\frac{1}{n}}}\\
&\qquad +O\Big(z^{r+1}\log z, z^{r+\frac{n-1}{n}}w(\log
z)^{\frac{n-1}{n}}\Big)\\
&\qquad  +O\Big( z^{r+\frac{n-3}{n}}w^{2}(\log
z)^{\frac{n-3}{n}},\ldots,z^rw^{n-1}\Big)\bigg]+z^{r+1+\frac{1}{n}}\psi_1(z),
\end{aligned}
\end{cases}
\end{equation}
for some $a\neq 0$.
\end{enumerate}
Furthermore, in both cases $\hat{f}_1$, $\hat{f}_2$ and $\psi_1(z)$
are analytic in $z$, $(\log z)^{\frac{1}{n}}$ and $w$, with $\psi_1$
of the form
\begin{displaymath}
\psi_1(z)=\sum_{j=1}^{\infty}z^{\frac{j-1}{n}}R_j^1\left((\log
z)^{\frac{1}{n}}\right),
\end{displaymath}
where $R_j^1(t)$ are analytic functions on the image, under the
transformation $t=(\log z)^{\frac{1}{n}}$, of our suitable simply
connected open set with $z=0$ on its boundary.
\end{prop}

\begin{proof}
Let $n\geq 2$. In the coordinates $(Z,W)$, with
\begin{equation}\label{1change}
\begin{cases}
\begin{aligned}
Z&=\alpha z,\\
W&=w,
\end{aligned}
\end{cases}
\end{equation}
and $\alpha\neq 0$ arbitrary parameter, \eqref{fstartn} becomes
\begin{equation}\label{barf}
\begin{cases}
\begin{aligned}
\bar{f}_1(Z,W)&=Z +a_{0,n-1}\alpha ^{-r} Z^{r+1}W^{n-1} +
a_{1,0}\alpha ^{-r-1}Z^{r+2} \\
&\qquad+O(Z^{r+3},Z^{r+2}W,Z^{r+1}W^{n}),\\
\bar{f}_2(Z,W)&=W[1 + b_{0,n}\alpha ^{-r}Z^rW^{n-1} +
b_{1,1}\alpha ^{-r-1}Z^{r+1}\\
&\qquad +O(Z^{r+2},Z^{r+1}W,Z^rW^n)] + b_{1,0}\alpha ^{-r-1}Z^{r+1}
+ O(Z^{r+2}).
\end{aligned}
\end{cases}
\end{equation}
If $n=2$, let us consider the following change of variables
\begin{equation}\label{2change}
\begin{cases}
u=Z,\\
v=W - aZ^{\frac{1}{2}}(\log Z)^{\frac{1}{2}},
\end{cases}
\end{equation}
with $a\neq 0$ arbitrary parameter. Then \eqref{barf} assumes the
following form
\begin{displaymath}
\begin{cases}
\begin{aligned}
\hat{f}_1(u,v)&=u +a_{0,1}\alpha ^{-r} u^{r+1}v + aa_{0,1}\alpha
^{-r}u^{r+1+\frac{1}{2}}(\log u)^{\frac{1}{2}} \\
&\quad \quad + O\Big(u^{r+2}\log u,u^{r+1+
\frac{1}{2}}v(\log u)^{\frac{1}{2}},u^{r+1}v^{2}\Big),\\
\hat{f}_2(u,v)&=v\bigg[1 + b_{0,2}\alpha ^{-r}u^rv + ab_{0,2}\alpha
^{-r}u^{r+\frac{1}{2}}(\log u)^{\frac{1}{2}} -\frac{aa_{0,1}\alpha
^{-r}}{2}\frac{u^{r+\frac{1}{2}}}{(\log u)^{\frac{1}{2}}} \\
&\quad \quad +O\Big(u^{r+1}\log u,u^{r+\frac{1}{2}}v(\log
u)^{\frac{1}{2}},u^rv^2\Big)\bigg]\\
&\quad \quad+ \alpha^{-r}\Big(
b_{1,0}\alpha^{-1}-a^2\frac{a_{0,1}}{2}\Big)u^{r+1}
+a^2\alpha^{-r}\Big(b_{0,2}-\frac{a_{0,1}}{2}\Big)u^{r+1}\log u\\
&\quad \quad+O\Big(u^{r+1+\frac{1}{2}}(\log u)^{1+\frac{1}{2}}\Big).
\end{aligned}
\end{cases}
\end{displaymath}
Since $a_{0,1}=2b_{0,2}$, if we choose $(\alpha,a)$ as solution of
\begin{displaymath}
\begin{cases}
\begin{aligned}
&b_{1,0}\alpha^{-1}-a^22^{-1}a_{0,1}=0,\\
&aa_{0,1}\alpha ^{-r}=-1,
\end{aligned}
\end{cases}
\end{displaymath}
we obtain \eqref{fhat}. It is not difficult to check that changes of
variables of type \eqref{1change} and \eqref{2change} generate in
$\hat{f}_2(u,0)$ terms of the types
\begin{displaymath}
\begin{cases}
\begin{aligned}
&u^{r+1+\frac{1}{2}}(\log u)^{1+\frac{1}{2}},\quad
u^{r+1+\frac{1}{2}}(\log u)^{\frac{1}{2}},\quad
\frac {u^{r+1+\frac{1}{2}}}{(\log u)^{\frac{1}{2}}},\\
&u^{r+2}(\log u)^2, \quad u^{r+2}\log u, \quad u^{r+2},\\
&u^{r+2+\frac{1}{2}}(\log u)^{1+\frac{1}{2}},\quad
u^{r+2+\frac{1}{2}}(\log u)^{\frac{1}{2}}, \quad
\frac {u^{r+2+\frac{1}{2}}}{(\log u)^{\frac{1}{2}}},\\
&u^{r+3}(\log u)^2, \quad u^{r+3}\log u, \quad u^{r+3},\ldots
\end{aligned}
\end{cases}
\end{displaymath}
then we can write
\begin{displaymath}
\begin{aligned}
\hat{f}_2(u,0)&=O\Big(u^{r+1+\frac{1}{2}}(\log
u)^{1+\frac{1}{2}}\Big)=u^{r+1+\frac{1}{2}}\Big[R^1_1\Big((\log
u)^{\frac{1}{2}}\Big)+u^{\frac{1}{2}}R^1_2\Big(
(\log u)^{\frac{1}{2}}\Big)\\
&\quad +uR^1_3\Big((\log
u)^{\frac{1}{2}}\Big)+u^{1+\frac{1}{2}}R^1_4\Big((\log
u)^{\frac{1}{2}}\Big)+\cdots\Big]=u^{r+1+\frac{1}{2}}\psi_1(u),
\end{aligned}
\end{displaymath}
where $R^1_j(t)$ are analytic functions on the image, under the
transformation $t=(\log u)^{\frac{1}{2}}$ of our suitable simply
connected open set with $u=0$ on its boundary.

Let $n\geq 3$. If we consider the following change of coordinates
\begin{displaymath}
\begin{cases}
u=Z,\\
v=W - aZ^{\frac{1}{n}}(\log Z)^{\frac{1}{n}},
\end{cases}
\end{displaymath}
with $a$ such that $a^n=(\alpha b_{0,n})^{-1}b_{1,0}$, then
\eqref{barf} becomes
\begin{displaymath}
\begin{cases}
\begin{aligned}
\hat {f}_1(u,v)&=u+(n-1)a^{n-2}a_{0,n-1}\alpha^{-r}
u^{r+1+\frac{n-2}{n}}v(\log u)^{\frac{n-2}{n}} \\
&\quad \,+a^{n-1}a_{0,n-1}\alpha^{-r}u^{r+1+\frac{n-1}{n}}(\log
u)^{\frac{n-1}{n}}\\
&\quad \,+O\Big(u^{r+2}\log u,u^{r+1+\frac{n-1}{n}}v(\log
u)^{\frac{n-1}{n}} \Big)\\
&\quad \, +O\Big(u^{r+1+\frac{n-3}{n}}v^{2}(\log
u)^{\frac{n-3}{n}},u^{r+1+\frac{n-4}{n}}v^{3}(\log
u)^{\frac{n-4}{n}},\ldots,u^{r+1}v^{n-1}\Big),\\
\hat{f}_2(u,v)&=v\Big[1+(n-1)a^{n-2}b_{0,n}\alpha^{-r}u^{r+
\frac{n-2}{n}}v(\log u)^{\frac{n-2}{n}}\\
&\quad \, +a^{n-1}b_{0,n}\alpha^{-r}u^{r+\frac{n-1}{n}}(\log
u)^{\frac{n-1}{n}} -\frac{n-1}{n}a^{n-1}a_{0,n-1}\alpha^{-r}
\frac{u^{r+\frac{n-1}{n}}}{(\log u)^{\frac{1}{n}}}\\
&\quad  \,+O\Big(u^{r+1}\log u,
u^{r+\frac{n-1}{n}}v(\log u)^{\frac{n-1}{n}}\Big) \\
&\quad  \,+O\Big(u^{r+\frac{n-3}{n}}v^{2}(\log
u)^{\frac{n-3}{n}},u^{r+\frac{n-4}{n}}v^{3}(\log
u)^{\frac{n-4}{n}},\ldots,u^rv^{n-1}\Big)
\Big]\\
&\quad \,+O\Big(u^{r+1+\frac{1}{n}}(\log u)^{1+\frac{1}{n}}\Big).
\end{aligned}
\end{cases}
\end{displaymath}
Since $a_{0,n-1}=nb_{0,n}$, choosing $(\alpha, a)$ as solution of
\begin{displaymath}
\begin{cases}
\begin{aligned}
&b_{1,0}\alpha^{-1}-a^nn^{-1}a_{0,n-1}=0,\\
&a^{n-1}a_{0,n-1}\alpha^{-r}=-1,
\end{aligned}
\end{cases}
\end{displaymath}
we obtain \eqref{fhatn}.

As before, analyzing pure terms in $\hat{f}_2(u,0)$ we have
\begin{displaymath}
\begin{aligned}
\hat{f}_2(u,0)&=O\left(u^{r+1+\frac{1}{n}}(\log
u)^{1+\frac{1}{n}}\right)=u^{r+1+\frac{1}{n}}\left[ R^1_1\left((\log
u)^{\frac{1}{n}}\right)+u^{\frac{1}{n}}R^1_2
\left((\log u)^{\frac{1}{n}}\right) \right. \\
&\left. \quad +u^{\frac{2}{n}}R^1_3\left((\log
u)^{\frac{1}{n}}\right) +\cdots +uR^1_{n+1}\left((\log
u)^{\frac{1}{n}}\right)+\cdots \right]=u^{r+1+\frac{1}{n}}\psi_1(u),
\end{aligned}
\end{displaymath}
where $R^1_j(t)$ are analytic functions on the image, under the
transformation $t=(\log u)^{\frac{1}{n}}$ of our suitable simply
connected open set with $u=0$ on its boundary.
\end{proof}

\begin{nota}
In particular, for any $j\geq 1$ we can write $R_j^1(t)$ in the
following form
\begin{displaymath}
R_j^1(t):=\sum_{k=1}^{\infty}\frac{c_{j,-k}^1}{t^k}+c^1_{j,0}+
c^1_{j,1}t+\cdots+c^1_{j,m^1_j}t^{m^1_j}=:S_j^1(t)+P_j^1(t),
\end{displaymath}
for suitable $c^1_{j,k}\in \C$ and with
$P_j^1(t):=c^1_{j,0}+c^1_{j,1}t+\cdots+c^1_{j,m^1_j}t^{m^1_j}$ such
that $\deg P_j^1=m^1_j$ or $P_j^1(t)\equiv 0$. We shall say that
$R_j^1$ satisfies the {\em star property}.
\end{nota}

\begin{prop}\label{shiftaggio}
Let $(\hat{f_1},\hat{f}_2)$ be a map of the form \eqref{fhat}, if
$n=2$, or \eqref{fhatn}, if $n\geq 3$. Then there exists a sequence
$\{Q_h(t)\}_{h\in \mathbb N}$ of functions analytic on the image,
under the transformation $t=(\log z)^{\frac{1}{n}}$, of our suitable
simply connected open set with $z=0$ on its boundary, such that if
$w_{h+1}(z)$ is defined by
\begin{displaymath}
\begin{cases}
w_1(z)=z^{\frac{2}{n}}Q_0\left((\log z)^{\frac{1}{n}}\right),\\
w_{h+1}(z)=w_h(z)+ z^{\frac{h+2}{n}}Q_h\left((\log
z)^{\frac{1}{n}}\right), \qquad \textrm{if $h\geq 1$},
\end{cases}
\end{displaymath}
then
\begin{equation}\label{h-shift-n}
\hat{f}_2(z,w_{h+1}(z))-w_{h+1}\left(\hat{f}_1(z,w_{h+1}(z))\right)
=z^{r+1+\frac{h+2}{n}}\psi_{h+2}(z),
\end{equation}
where $\psi_{h+2}$ is of the form
$\psi_{h+2}(z)=\sum_{j=1}^{\infty}z^{\frac{j-1}{n}}R_j^{h+2}\left((\log
z )^{\frac{1}{n}}\right)$, with $R_j^{h+2}(t)$ analytic functions on
the image, under the transformation $t=(\log z )^{\frac{1}{n}}$, of
our suitable simply connected open set with $z=0$ on its boundary.
\end{prop}

\begin{proof}
In both cases we know that
\begin{displaymath}
\psi_1(z)=R_1^1\left((\log
z)^{\frac{1}{n}}\right)+z^{\frac{1}{n}}R_2^1\left((\log
z)^{\frac{1}{n}}\right)+\cdots+zR_{n+1}^1\left((\log
z)^{\frac{1}{n}} \right)+\cdots,
\end{displaymath}
where $R_j^1(t)$ are analytic functions on the image, under the
transformation $t=(\log z)^{\frac{1}{n}}$ of our suitable simply
connected open set with $z=0$ on its boundary.

We prove the proposition by induction on $h$.\\
Let $h=0$. If we define $w_1(z):=z^{\frac{2}{n}}Q_{0}\left((\log
z)^{\frac{1}{n}}\right)$, with $Q_0(t)$ a holomorphic solution of
the differential equation
\begin{equation}\label{1eqdiff}
t^{-(n-1)}Q'(t)+ (n-1)\left(1+\frac {1}{t^n}\right)Q(t)=
-nt^{-(n-1)}R_1^1(t),
\end{equation}
we have
\begin{displaymath}
\hat{f}_2(z,w_1(z))-w_1\left(\hat{f}_1(z,w_1(z))\right)=
z^{r+1+\frac{2}{n}}\psi_2(z),
\end{displaymath}
with
\begin{displaymath}
\psi_2(z)=R_1^2\left((\log z)^{\frac{1}{n}}\right)+
z^{\frac{1}{n}}R_2^2\left((\log z)^{\frac{1}{n}}\right)+\cdots
+zR_{n+1}^2\left((\log z)^{\frac{1}{n}}\right)+\cdots ,
\end{displaymath}
and where $R_j^2(t)$ are analytic functions on the image, under the
transformation $t=(\log z)^{\frac{1}{n}}$, of our suitable simply
connected open set with $z=0$ on its boundary.

We study the differential equation \eqref{1eqdiff}.

\noindent First, an easy computation shows that there exists a
unique formal solution of the following form
\begin{displaymath}
F_0(t):=\sum_{k=1}^{\infty}\frac{d^0_{-k}}{t^k}+d^0_0+d^0_1t
+\cdots+d^0_{l_0}t^{l_0},
\end{displaymath}
with $T^0(t):=d^0_0+d^0_1t+\cdots+d^0_{l_0}t^{l_0}$ such that $\deg
T^0=l_0$ or $T^0(t)\equiv 0$.

\noindent According to \cite{Wa} (Theorem 12.1, p. 57) there exists,
for sufficiently large $t$ in an open sector $S$ with vertex at the
origin and a positive central angle not exceeding $\frac{\pi}{n}$, a
holomorphic solution $Q_0(t)$ of \eqref{1eqdiff}, which is
asymptotic to the formal solution $F_0(t)$, in every proper
subsector $S'$ of $S$. This means that for all $m\geq 1$
\begin{equation}\label{asymp}
\lim_{t\to \infty }
t^m\bigg[Q_0(t)-\bigg(\sum_{k=1}^{m}\frac{d^0_{-k}}{t^k}+d^0_0+d^0_1t
+\cdots+d^0_{l_0}t^{l_0}\bigg) \bigg]=0,
\end{equation}
as $t$ tends to $\infty$ in $S'$, and we write $Q_0(t)\sim F_0(t)$
for short.

Then \eqref{h-shift-n} is satisfied for $h=0$.

Moreover, using Theorems 8.2, 8.3 and 8.8 of \cite{Wa} it follows
that $R^2_1(t)$ is asymptotic, in every proper subsector $S'$ of
$S$, to
\begin{displaymath}
\sum_{k=1}^{\infty}\frac{c_{1,-k}^2}{t^k}+c^2_{1,0}+c^2_{1,1}t+
\cdots+c^2_{1,m^2_1}t^{m^2_1}=:S_1^2(t)+P_1^2(t),
\end{displaymath}
for suitable $c^2_{1,k}\in \C$ and with
$P_1^2(t):=c^2_{1,0}+c^2_{1,1}t+\cdots+c^2_{1,m^2_1}t^{m^2_1}$ such
that $\deg P_1^2=m^2_1$ or $P_1^2(t)\equiv 0$.

For $h>0$, one sees by induction that if \eqref{h-shift-n} is
satisfied by $w_{h}$, \eqref{h-shift-n} is then satisfied by
$w_{h+1}(z)=w_{h}(z)+z^{\frac{h+2}{n}}Q_{h}\left((\log
z)^{\frac{1}{n}}\right)$ for a function $Q_h$, if and only if $Q_h$
is a solution of the differential equation
\begin{equation}\label{eqdiffconh}
t^{-(n-1)}Q'(t)+(n-1)\left[ \left(
h+1\right)+\frac{1}{t^n}\right]Q(t)=-nt^{-(n-1)}R_1^{h+1}(t).
\end{equation}
\noindent

We observe that if $R_1^{h}$ and $Q_{h-1}$ are asymptotic to some
functions which satisfy the star property, then also $R_1^{h+1}$ is
asymptotic to a function which satisfies the same property, for all
$h\geq 1$. Then, by the inductive hypothesis we can suppose that
$R_1^{h+1}(t)$ is asymptotic to one function which satisfies the
star property, hence
\begin{equation}\label{R^h+1}
R_1^{h+1}(t)\sim\sum_{k=1}^{\infty}\frac{c_{1,-k}^{h+1}}{t^k}+
c^{h+1}_{1,0}+c^{h+1}_{1,1}t+\cdots+c^{h+1}_{1,m^{h+1}_1}t^{m^{h+1}_1}
=:S_1^{h+1}(t)+P_1^{h+1}(t),
\end{equation}
for suitable $c^{h+1}_{1,k}\in \C$ and with
$P_1^{h+1}(t):=c^{h+1}_{1,0}+c^{h+1}_{1,1}t+
\cdots+c^{h+1}_{1,m^{h+1}_1}t^{m^{h+1}_1}$ such that $\deg
P_1^{h+1}=m^{h+1}_1$ or $P_1^{h+1}(t)\equiv 0$. Also in this case it
is possible to prove that there exists a unique formal solution of
the differential equation \eqref{eqdiffconh} of the form
\begin{displaymath}
F_h(t):=\sum_{k=1}^{\infty}\frac{d^{h}_{-k}}{t^k}+d^h_0+d^h_1t+
\cdots+d^h_{l_h}t^{l_h},
\end{displaymath}
with $T^h(t):=d^h_0+d^h_1t+\cdots+d^h_{l_h}t^{l_h}$ such that $\deg
T^h=l_h$ or $T^h(t)\equiv 0$.

According to \cite{Wa} (Theorem 12.1, p. 57) there exists, for
sufficiently large $t$ in a suitable open sector $S$ with vertex at
the origin and a positive central angle not exceeding
$\frac{\pi}{n}$, a holomorphic solution $Q_h(t)$ of
\eqref{eqdiffconh} which is asymptotic to the formal solution
$F_h(t)$, in every proper subsector $S'$ of $S$. In particular, this
means that for all $m\geq 1$
\begin{displaymath}
Q_h(t)=\sum_{k=1}^{m}\frac{d^h_{-k}}{t^k}+d^h_0+d^h_1t+\cdots
+d^h_{l_h}t^{l_h}+O\Big(\frac{1}{t^m}\Big),
\end{displaymath}
as $t$ tends to $\infty$ in $S'$.
\end{proof}

\begin{corollario}\label{shift}
Let $(\hat{f_1},\hat{f}_2)$ be a map of the form \eqref{fhat}, if
$n=2$, or \eqref{fhatn}, if $n\geq 3$. Then for every $h\in\N$, we
can choose local coordinates relative to which the map takes the
form:
\begin{enumerate}
\item if $n=2$:
\begin{equation}\label{xshift}
\begin{cases}
\begin{aligned}
\hat{f}_1(Z,W)&=Z -a^{-1} Z^{r+1}W -Z^{r+1+\frac{1}{2}}(\log
Z)^{\frac{1}{2}}\\
&\qquad + O\Big(Z^{r+2}(\log Z)^i,Z^{r+1+\frac{1}{2}}W(\log
Z)^{\frac{1}{2}},Z^{r+1}W^{2}\Big),\\
\hat{f}_2(Z,W)&=W\bigg[1 -\frac{a^{-1}}{2}Z^rW -\frac{1}{2}
Z^{r+\frac{1}{2}}(\log Z)^{\frac{1}{2}} +\frac{1}{2}
\frac{Z^{r+\frac{1}{2}}}{(\log Z)^{\frac{1}{2}}} \\
&  \qquad  + O\left(Z^{r+1}(\log Z)^i,Z^{r+\frac{1}{2}}W(\log
Z)^{\frac{1}{2}},Z^rW^2\right)\bigg]\\
& \qquad  + Z^{r+1+\frac{h+2}{2}}\psi_{h+2}(Z),
\end{aligned}
\end{cases}
\end{equation}
for some $i\geq 1$, where $\hat{f}_1$, $\hat{f}_2$ and
$\psi_{h+2}(Z)$ are analytic in $Z$, $(\log Z)^{\frac{1}{2}}$, and
$W$, with $Z$ belonging to our suitable simply connected open set
with $Z=0$ on the boundary; \item if $n\geq 3$:
\begin{equation}\label{xshiftn}
\begin{cases}
\begin{aligned}
\hat {f}_1(Z,W)&=Z-(n-1)a^{-1}Z^{r+1+\frac{n-2}{n}}W(\log
Z)^{\frac{n-2}{n}} - Z^{r+1+\frac{n-1}{n}}(\log
Z)^{\frac{n-1}{n}}\\
&\qquad+ O\left(Z^{r+2}(\log Z)^i,Z^{r+1+\frac{n-1}{n}}W(\log
Z)^{\frac{n-1}{n}} \right)\\
&\qquad +O\left(Z^{r+1+\frac{n-3}{n}}W^{2}(\log
Z)^{\frac{n-3}{n}},\ldots,Z^{r+1}W^{n-1}\right),\\
\hat{f}_2(Z,W)&=W\bigg[1-\frac{n-1}{n}a^{-1}Z^{r+\frac{n-2}{n}}
W(\log Z)^{\frac{n-2}{n}} -\frac{1}{n}Z^{r+\frac{n-1}{n}}
(\log Z)^{\frac{n-1}{n}}\\
&\qquad +\frac{n-1}{n}\frac{Z^{r+\frac{n-1}{n}}}{(\log
Z)^{\frac{1}{n}}}+O\left(Z^{r+1}(\log Z)^i,Z^{r+\frac{n-1}{n}}
W(\log Z)^{\frac{n-1}{n}}\right)\\
&\qquad +O\left(Z^{r+\frac{n-3}{n}}W^{2}(\log
Z)^{\frac{n-3}{n}},\ldots,Z^rW^{n-1}\right)
\bigg]\\
&\qquad  +Z^{r+1+\frac{h+2}{n}}\psi_{h+2}(Z),
\end{aligned}
\end{cases}
\end{equation}
for some $i\geq 1$, where $\hat{f}_1$, $\hat{f}_2$ and
$\psi_{h+2}(Z)$ are analytic in $Z$, $(\log Z)^{\frac{1}{n}}$, and
$W$, with $Z$ belonging to our suitable simply connected open set
with $Z=0$ on the boundary.
\end{enumerate}
\end{corollario}

\begin{proof}
Let $w_{h+1}$ be given by Proposition \ref{shiftaggio}. One has just
to take as new coordinates
\begin{displaymath}
\begin{cases}
Z=z,\\
W=w-w_{h+1}(z).
\end{cases}
\end{displaymath}
We remark that $i=\max \{1,I/n\}$, where $I$ is the maximum power of
$t$ appearing in the following estimate, which is obtained from
\eqref{asymp} and is satisfied for all $m\geq 1$ and for
sufficiently large $t$,
\begin{displaymath}
Q_0(t)=\sum_{k=1}^m\frac{d^0_{-k}}{t^k}+d^0_0+d^0_1t+\cdots+
d^0_{l_0}t^{l_0}+O\Big(\frac{1}{t^m}\Big).
\end{displaymath}
Hence $I=l_0$ or $I=-\min \{k\geq 1 \mid d^0_{-k}\neq 0\}$.
\end{proof}

By Corollary \ref{shift}, we can choose local coordinates defined in
our suitable simply connected open set with $O$ on the boundary,
such that for every $h\in\N$, the map $\hat{f}$ is written as in
\eqref{xshift} if $n=2$ or \eqref{xshiftn} if $n\geq 3$. Let us
consider $h=2n-3$; for our purposes it will be enough.

Now set
\begin{displaymath}
D_{r+\frac{n-1}{n},\delta}:=\left\{\zeta\in \mathbb C\quad |\quad
|\zeta^{r+\frac{n-1}{n}} (\log \zeta)^{\frac{n-1}{n}} - \delta |<
\delta \right\},
\end{displaymath}
for some $\delta>0$.

For an arbitrary $\delta$, this set has a number of connected
components which depends on $r$ and $n$: if $r$ is even,
$D_{r+\frac{n-1}{n},\delta}$ has $r+1$ connected and simply
connected components with the origin on the boundary; if $r$ is odd
but $n$ is not too big, the number of the connected components with
the origin on the boundary is $r$, otherwise if $n$ is big enough
this number is $r+1$.

Instead, if $\delta$ is small enough, the number of the components
with the origin on the boundary is equal to $r+1$, in all cases, and
the picture of the set looks like a flower, in which not all the
petals have the same size.

On these components, and for $\delta$ sufficiently small, the
functions $w_h$ are well defined.

Let us fix some notation. Put
\begin{displaymath}
\mathcal{E}_{n}^{\beta}(\delta):= \left\{ w\in \mbox{
Hol}\big(D_{r+\frac{n-1}{n},\delta} , \mathbb C\big) \quad | \quad
w(\zeta)=\zeta^{2}(\log \zeta)^{\beta} h^o(\zeta), \quad
\|h^o\|_\infty < \infty  \right\},
\end{displaymath}
with $\beta \in \R$. It is a Banach space with the norm
$\|w\|_{\mathcal{E}_{n}^{\beta}(\delta)}=\|h^o\|_{\infty}$. For
$w\in\mathcal{E}_{n}^{\beta}(\delta)$, put
$\hat{f}^w(\zeta)=\hat{f}_1(\zeta,w(\zeta))$.

Assume we have found $\hat{w}\in\mathcal{E}_{n}^{\beta}(\delta)$
which satisfies the following properties:

\begin{enumerate}
\item[(i)] for all $\zeta \in D_{r+\frac{n-1}{n},\delta}$
\begin{equation}\label{eqfunz}
\hat{w}\left(\hat{f}_1(\zeta,\hat{w}(\zeta))\right)=
\hat{f}_2(\zeta,\hat{w}(\zeta));
\end{equation}

\item[(ii)] there exists a positive constant $\delta_0$ such that,
if $0<\delta<\delta_0$ then $\hat{f}^{\hat{w}}$ sends every
component of $D_{r+\frac{n-1}{n},\delta}$ into itself and
$\hat{\zeta}_k\to 0$ as $k\to \infty$, where $\hat{\zeta}_k$ denotes
the sequence of iterate of an arbitrary point $\zeta\in
D_{r+\frac{n-1}{n},\delta}$ by the
transformation$\hat{f}^{\hat{w}}(\zeta):=\hat{\zeta}_1$;
\end{enumerate}
then the restriction of
$\hat{\varphi}(\zeta):=(\zeta,\hat{w}(\zeta))$ to any component of
$D_{r+\frac{n-1}{n},\delta}$ is a parabolic curve for $\hat{f}$. In
fact, equation \eqref{eqfunz} implies that
$\hat{\varphi}(D_{r+\frac{n-1}{n},\delta})$ is invariant under
$\hat{f}$ and property (ii) that
$\big(\hat{f}|_{\hat{\varphi}(D_{r+\frac{n-1}{n},\delta})}\big)^k\to
O$, as $k\to\infty$.

First, we must find a solution of \eqref{eqfunz}. The function
$\hat{w}$ will be obtained as fixed point of a functional operator
$T$, which acts as a contraction on a suitable closed convex subset
of $\mathcal{E}_{n}^{\beta}(\delta)$ and which we are going to
describe.

If $\hat{f}$ is given by \eqref{xshift} or \eqref{xshiftn} we denote
$z_1:=\hat{f}_1(z,w)$ and $w_1:=\hat{f}_2(z,w)$. Suppose $z$, $z_1$
belong to the same connected component of
$D_{r+\frac{n-1}{n},\delta}$ and define
\begin{displaymath}
H(z,w)=w-\frac{z^{\frac{1}{n}}}{z_1^{\frac{1}{n}}}w_1.
\end{displaymath}
A direct computation shows that:
\begin{enumerate}
\item if $n=2$
\begin{equation}\label{OH1/2}
H(z,w)=O\left( \frac{z^{r+\frac{1}{2}}}{(\log z)^{\frac{1}{2}}}w,
z^{r+\frac{1}{2}}w^2(\log z)^{\frac{1}{2}}, z^{r+2+\frac{1}{2}}(\log
z)^J \right),
\end{equation}
\item if $n\geq 3$
\begin{equation}\label{OHn}
\begin{aligned}
H(z,w)&=O\bigg( \frac{z^{r+\frac{n-1}{n}}}{(\log z)^{\frac{1}{n}}}w,
z^{r+\frac{n-1}{n}}w^2(\log z)^{\frac{n-1}{n}},
z^{r+2+\frac{n-1}{n}}(\log z)^J \bigg)\\
&\,\, +O\Big(z^{r+\frac{n-3}{n}}w^3(\log
z)^{\frac{n-3}{n}},z^{r+\frac{n-4}{n}}w^4(\log
z)^{\frac{n-4}{n}},\ldots,z^rw^n\Big),
\end{aligned}
\end{equation}
\end{enumerate}
where, in both cases, $J$ is the maximum power of $\log z$ appearing
in \eqref{R^h+1}, the asymptotic expansion of $R_1^{2n-1}\left((\log
z)^{\frac{1}{n}}\right)$.

\begin{nota}
In particular, if $\deg P^{2n-1}_1=m^{2n-1}_1$ then
$J=m^{2n-1}_1/n$, otherwise if $P^{2n-1}_1\equiv 0$ then $J=-\min
\{k\geq 1 \mid c^{2n-1}_{1,-k}\neq 0 \}/n$.
\end{nota}

Let us take $\beta=|J|$, and let $T$ be the operator on
$\mathcal{E}_{n}^{|J|}(\delta)$ defined by
\begin{equation}\label{Tn}
Tw(\zeta_0)=\zeta_0^{\frac{1}{n}}\sum_{k=0}^{\infty}
\zeta_k^{-\frac{1}{n}}H(\zeta_k,w(\zeta_k)),
\end{equation}
where $\zeta_k:=(\hat{f}^{w})^k(\zeta_0)$. In the remainder of this
section, we shall prove that $T$ is well defined for any $\zeta_0\in
D_{r+\frac{n-1}{n},\delta}$ and for $\delta$ sufficiently small, and
that restricted to a suitable closed convex subset of
$\mathcal{E}_{n}^{|J|}(\delta)$ is contracting, so that $T$ has a
fixed point. Assume that this is done, and let $\tilde{w}$ be the
fixed point. Then $\tilde{w}=\hat{w}$, since it is not difficult to
see that $\tilde{w}$ is a fixed point of $T$ if and only if it
satisfies \eqref{eqfunz}. We show the implication ($\Rightarrow$):
let $\zeta_0\in D_{r+\frac{n-1}{n},\delta}$, then
\begin{displaymath}
\begin{aligned}
\tilde{w}(\zeta_1)&=T\tilde{w}(\zeta_1)=
\zeta_1^{\frac{1}{n}}\sum_{k=1}^{\infty}\zeta_k^{-\frac{1}{n}}
H(\zeta_k,\tilde{w}(\zeta_k))\\
&=\frac{\zeta_1^{\frac{1}{n}}}{\zeta_0^{\frac{1}{n}}}\left
(\zeta_0^{\frac{1}{n}}\sum_{k=0}^{\infty}\zeta_k^{-\frac{1}{n}}
H(\zeta_k,\tilde{w}(\zeta_k))-H(\zeta_0,\tilde{w}(\zeta_0))\right)\\
&=\frac{\zeta_1^{\frac{1}{n}}}{\zeta_0^{\frac{1}{n}}}\Big
[T\tilde{w}(\zeta_0)-H(\zeta_0,\tilde{w}(\zeta_0))\Big ]=
\frac{\zeta_1^{\frac{1}{n}}}{\zeta_0^{\frac{1}{n}}}\left
[\tilde{w}(\zeta_0)-\tilde{w}(\zeta_0) +
\frac{\zeta_0^{\frac{1}{n}}}{\zeta_1^{\frac{1}{n}}}
\hat{f}_2(\zeta_0,\tilde{w}(\zeta_0))\right]\\
&=\hat{f}_2(\zeta_0,\tilde{w}(\zeta_0)).
\end{aligned}
\end{displaymath}

The next results establish some useful bounds. In particular, the
next lemma shows that property (ii) is satisfied for each
$w\in\mathcal{E}_{n}^{|J|}(\delta)$, if $\delta$ is chosen small
enough and if $\|h^o\|_{\infty}\leq 1$.

\begin{lemma}\label{ineq}
Let $\hat{f}$ be a map of the form \eqref{xshift} or
\eqref{xshiftn}, with $h=2n-3$. For each
$w\in\mathcal{E}_{n}^{|J|}(\delta)$ such that $\|h^o\|_{\infty}\leq
1$, let $\{\zeta_k\}$ be the sequence of iterates of $\zeta$ by the
transformation $\hat{f}^w(\zeta)=\hat{f}_1(\zeta,w(\zeta))$. Then
there is a positive constant $\delta_0$ such that if $\zeta$ belongs
to a connected component of $D_{r+\frac{n-1}{n},\delta_0}$ then
$\zeta_k$ belongs to the same component of $D_{r+\frac{n-1}
{n},\delta_0}$ for every $k\in\N$, and
\begin{equation}\label{inequalities}
\begin{aligned}
\frac{2}{3}\frac{|\zeta|^{r+\frac{n-1}{n}}|\log
\zeta|^{\frac{n-1}{n}}}{\big|1+k\left(r+\frac{n-1}{n}\right)
\zeta^{r+\frac{n-1}{n}}(\log \zeta)^{\frac{n-1}{n}}\big|} &\leq
|\zeta_k|^{r+\frac{n-1}{n}}|\log \zeta_k|^{\frac{n-1}{n}}\\
&\leq 2\frac{|\zeta|^{r+\frac{n-1}{n}}|\log
\zeta|^{\frac{n-1}{n}}}{\big|1+k\left(r+\frac{n-1}{n}\right)
\zeta^{r+\frac{n-1}{n}}(\log \zeta)^{\frac{n-1}{n}}\big|}.
\end{aligned}
\end{equation}
\end{lemma}

\begin{proof}
From the hypothesis on $w(\zeta)$, we have
\begin{displaymath}
\frac{1}{\zeta_1^{r+\frac{n-1}{n}}(\log \zeta_1)^{\frac{n-1}{n}}}=
\frac{1}{\zeta^{r+\frac{n-1}{n}}(\log
\zeta)^{\frac{n-1}{n}}}+r+\frac{n-1}{n}+g(\zeta)
\end{displaymath}
where $g(\zeta)=O\left(\frac{1}{\log \zeta}\right)=o(1)$. Observe
that $\zeta$ and $\zeta_1$ are in the same component, i.e. the
desired component is invariant under $\hat{f}$. Moreover, $g(\zeta)$
represents a function bounded by $K\frac{1}{|\log \zeta|}$, where
$K$ is a constant independent of $w$. Iterating and dividing by $k$,
we get
\begin{equation}\label{rif}
\frac{1}{k\zeta_k^{r+\frac{n-1}{n}}(\log \zeta_k)^{\frac{n-1}{n}}}=
\frac{1}{k\zeta^{r+\frac{n-1}{n}}(\log
\zeta)^{\frac{n-1}{n}}}+\left(r+\frac{n-1}{n}\right)+
\frac{\sum_{j=0}^{k-1}g(\zeta_j)}{k},
\end{equation}
where $g_k(\zeta)=\frac{1}{k}\sum_{j=0}^{k-1}g(\zeta_j)$ converges
to $0$ when $k$ goes to $\infty$, by Ces\`aro convergence. So, with
$\zeta\in D_{r+\frac{n-1}{n},\delta}$ and $\delta$ sufficiently
small, we have
\begin{equation}\label{andasint}
\left(r+\frac{n-1}{n}\right)\zeta_k^{r+\frac{n-1}{n}} (\log
\zeta_k)^{\frac{n-1}{n}}\sim \frac{1}{k}.
\end{equation}

From \eqref{rif} we have, $\forall k\in \N$
\begin{equation}\label{prin}
\begin{aligned}
\frac{1}{\zeta_k^{r+\frac{n-1}{n}}(\log \zeta_k)^{\frac{n-1}{n}}}&=
\frac{1}{\zeta^{r+\frac{n-1}{n}}(\log
\zeta)^{\frac{n-1}{n}}}\left(1+k\left(r+\frac{n-1}{n}\right)
\zeta^{r+\frac{n-1}{n}}(\log \zeta)^{\frac{n-1}{n}}\right)\\
&\quad \cdot \left[1+ \frac{\zeta^{r+\frac{n-1}{n}}(\log
\zeta)^{\frac{n-1}{n}}\sum_{j=0}^{k-1}g(\zeta_j)}
{1+k\left(r+\frac{n-1}{n}\right)\zeta^{r+\frac{n-1}{n}}(\log
\zeta)^{\frac{n-1}{n}}}\right],
\end{aligned}
\end{equation}
where $\sum_{j=0}^{k-1}g(\zeta_j)=o(k)$. Hence, for $\zeta\to 0$
\begin{displaymath}
\frac{\zeta^{r+\frac{n-1}{n}}(\log
\zeta)^{\frac{n-1}{n}}\sum_{j=0}^{k-1}g(\zeta_j)}
{1+k\left(r+\frac{n-1}{n}\right)\zeta^{r+\frac{n-1}{n}}(\log
\zeta)^{\frac{n-1}{n}}} \to 0
\end{displaymath}
and the quantities inside the brackets in \eqref{prin} are uniformly
close to $1$; so for $\delta$ small enough we get the inequalities
\eqref{inequalities}.
\end{proof}

\begin{corollario}\label{corserie}
Let $\{\zeta_k\}$ and $\delta=\delta_0$ be defined as in Lemma
\ref{ineq}, with $w\in\mathcal{E}_{n}^{|J|}(\delta)$ such that
$\|h^o\|_{\infty}\leq 1$. Then for each real number
$s>r+\frac{n-1}{n}$, for each $q\geq 0$, there exists a constant
$C_{s,q}$ such that, for every $\zeta\in
D_{r+\frac{n-1}{n},\delta}$, we have
\begin{displaymath}
\sum_{k=0}^{\infty}|\zeta_k|^s|\log \zeta_k|^q\leq C_{s,q}
|\zeta|^{s-(r+\frac{n-1}{n})}|\log |\zeta||^{q-\frac{n-1}{n}}.
\end{displaymath}
\end{corollario}

\begin{proof}
Let $\zeta\in D_{r+\frac{n-1}{n},\delta}$. From the first inequality
of \eqref{inequalities} we get, for all $k\in \N$
\begin{displaymath}
\frac{2}{3}\frac{|\zeta|^{r+\frac{n-1}{n}}}{\big|1+
k\left(r+\frac{n-1}{n}\right)\zeta^{r+\frac{n-1}{n}}(\log
\zeta)^{\frac{n-1}{n}}\big|}\frac{1}{|\log \zeta_k|^{\frac{n-1}{n}}}
\leq |\zeta_k|^{r+\frac{n-1}{n}},
\end{displaymath}
and for $\delta$ small enough
\begin{displaymath}
\log \left( \frac{2|\zeta|^{r+\frac{n-1}{n}}|\log
\zeta_k|^{\frac{1-n}{n}}}{3\big|1+k\left(r+
\frac{n-1}{n}\right)\zeta^{r+\frac{n-1}{n}}(\log
\zeta)^{\frac{n-1}{n}}\big|} \right) \leq
\left(r+\frac{n-1}{n}\right)\log |\zeta_k|<0.
\end{displaymath}
Taking the modulus, we get
\begin{displaymath}
\begin{aligned}
\left(r+\frac{n-1}{n}\right)&|\log |\zeta_k||\\
&\leq \Bigg| \log \left
(\frac{2}{3}\frac{|\zeta|^{r+\frac{n-1}{n}}}{\big|1+
k\left(r+\frac{n-1}{n}\right)\zeta^{r+\frac{n-1}{n}} (\log
\zeta)^{\frac{n-1}{n}}\big|}\frac{1}
{|\log \zeta_k|^{\frac{n-1}{n}}} \right) \Bigg|\\
&\leq \Bigg| \log \left(\frac{2}{3}\frac{|\zeta|^{r+\frac{n-1}{n}}}
{\big|1+k\left(r+\frac{n-1}{n}\right)\zeta^{r+\frac{n-1}{n}} (\log
\zeta)^{\frac{n-1}{n}}\big|} \right)
\Bigg|\\
&\quad +\frac{n-1}{n}|\log|\log \zeta_k||.
\end{aligned}
\end{displaymath}
Since $\log|\log z|\leq|\log|z||$, for $z$ sufficiently close to
$0$, the last inequality becomes
\begin{equation}\label{ineqlog}
\begin{aligned}
r|\log |\zeta_k||&\leq \Bigg| \log \left(
\frac{2}{3}\frac{|\zeta|^{r+\frac{n-1}{n}}}{\big|1+
k\left(r+\frac{n-1}{n}\right)\zeta^{r+\frac{n-1}{n}}(\log
\zeta)^{\frac{n-1}{n}}\big|} \right) \Bigg|\\
&\leq \left(r+\frac{n-1}{n}\right)|\log |\zeta||+\Bigg| \log
\frac{\big|1+k\left(r+\frac{n-1}{n}\right)\zeta^{r+\frac{n-1}{n}}
(\log \zeta)^{\frac{n-1}{n}}\big|}{\frac{2}{3}} \Bigg|.
\end{aligned}
\end{equation}
Since $\zeta\in D_{r+\frac{n-1}{n},\delta}$ then
$\hbox{Re}\left(\zeta^{r+\frac{n-1}{n}}(\log
\zeta)^{\frac{n-1}{n}}\right)>0$, and we get
\begin{displaymath}
\begin{aligned}
0<\log\frac{\big|1+k\left(r+\frac{n-1}{n}\right)
\zeta^{r+\frac{n-1}{n}}(\log \zeta)^{\frac{n-1}{n}}\big|}
{\frac{2}{3}}\\
\leq\log\frac{1+\big|k\left(r+\frac{n-1}{n}\right)
\zeta^{r+\frac{n-1}{n}}(\log
\zeta)^{\frac{n-1}{n}}\big|}{\frac{2}{3}}.
\end{aligned}
\end{displaymath}
Therefore, we can write \eqref{ineqlog} as
\begin{displaymath}
|\log |\zeta_k||\leq \frac{1}{r}\left[\left(r+
\frac{n-1}{n}\right)|\log
|\zeta||+\log\frac{1+\big|k\left(r+\frac{n-1}{n}\right)
\zeta^{r+\frac{n-1}{n}}(\log
\zeta)^{\frac{n-1}{n}}\big|}{\frac{2}{3}} \right].
\end{displaymath}

From the second inequality of \eqref{inequalities} we have, for all
$k\in \N$
\begin{displaymath}
|\zeta_k|^{r+\frac{n-1}{n}}\leq
2\frac{|\zeta|^{r+\frac{n-1}{n}}}{\big|1+
k\left(r+\frac{n-1}{n}\right)\zeta^{r+\frac{n-1}{n}}(\log
\zeta)^{\frac{n-1}{n}}\big|}\frac{|\log \zeta|^{\frac{n-1}{n}}}
{|\log \zeta_k|^{\frac{n-1}{n}}}.
\end{displaymath}
From \eqref{andasint} it follows that there exists a constant $C$
independent of $k$ and $\zeta$, such that $|\log \zeta|\leq C|\log
\zeta_k|$, for all $k\in\N$. Hence we have, for all $s>0$ and
$k\in\N$
\begin{displaymath}
|\zeta_k|^s\leq C_s\frac{|\zeta|^s}{\big|1+k\left(r+
\frac{n-1}{n}\right)\zeta^{r+\frac{n-1}{n}}(\log
\zeta)^{\frac{n-1}{n}}\big|^{s(r+\frac{n-1}{n})^{-1}}},
\end{displaymath}
for some constant $C_s$ depending only on $s$. Since
$\hbox{Re}\left(\zeta^{r+\frac{n-1}{n}}(\log
\zeta)^{\frac{n-1}{n}}\right)>0$, we get
\begin{displaymath}
|\zeta_k|^s\leq
C_s\frac{|\zeta|^s}{\left(1+\big|k\left(r+\frac{n-1}{n}\right)
\zeta^{r+\frac{n-1}{n}}(\log
\zeta)^{\frac{n-1}{n}}\big|^2\right)^{\frac{s}{2}
(r+\frac{n-1}{n})^{-1}}}.
\end{displaymath}
Hence, there are constants $K_1$, $K_{s,q}$ and $K_{s,q}^1$ such
that
\begin{displaymath}
\begin{aligned}
&\sum_{k=0}^{\infty}|\zeta_k|^s|\log \zeta_k|^q\leq K_1
\sum_{k=0}^{\infty}|\zeta_k|^s|\log |\zeta_k||^q\\
&\leq K_{s,q}|\zeta|^s \sum_{k=0}^{\infty} \frac{\left[
\left(r+\frac{n-1}{n}\right)|\log |\zeta||+ \log
\frac{1+\big|k\left(r+\frac{n-1}{n}\right)\zeta^{r+\frac{n-1}{n}}
(\log \zeta)^{\frac{n-1}{n}}\big|}{\frac{2}{3}}\right]^q}
{\left(1+\big|k\left(r+\frac{n-1}{n}\right)\zeta^{r+\frac{n-1}{n}}
(\log \zeta)^{\frac{n-1}{n}}\big|^{2}\right)^{\frac{s}{2}
(r+\frac{n-1}{n})^{-1}}}\\
&\leq K_{s,q}^{1}|\zeta|^{s-(r+\frac{n-1}{n})}|\log
\zeta|^{-\frac{n-1}{n}}
\int_{0}^{\infty}\frac{\left[\left(r+\frac{n-1}{n}\right)|\log
|\zeta||+\log \frac{3(1+u)}{2}\right]^q}{(1+u^2)^{\frac{s}{2}
(r+\frac{n-1}{n})^{-1}}}du.
\end{aligned}
\end{displaymath}
The last integral converges if and only if $s>r+\frac{n-1}{n}$ and
it gives us a contribute of type $K_{s,q}^2|\log \zeta|^q$, for some
constant $K_{s,q}^2$.
\end{proof}

Let us now study more carefully the operator $T$. From
\eqref{OH1/2}, or equivalently \eqref{OHn}, and Corollary
\ref{corserie} we get the following estimate
\begin{displaymath}
|Tw(\zeta)|\leq C_{n,J}|\zeta|^2|\log|\zeta||^{|J|-\frac{n-1}{n}}<
|\zeta|^2|\log|\zeta||^{|J|},
\end{displaymath}
hence $T$ sends the unit ball of $\mathcal{E}_{n}^{|J|}(\delta)$
into itself. We will be able to prove more than this: $T$ sends into
itself the convex closed set
\begin{displaymath}
\mathcal{F}_{n}^{|J|}(\delta)=\left\{w\in\mathcal{E}_{n}^{|J|}
(\delta)\quad|\quad |w(\zeta)|\leq|\zeta|^2|\log|\zeta||^{|J|},\quad
|w'(\zeta)|\leq|\zeta||\log |\zeta||^{|J|}\right\}.
\end{displaymath}
Then it will be enough to show that $T$ is a contraction on
$\mathcal{F}_{n}^{|J|}(\delta)$. But first we need another set of
results.

\begin{lemma}\label{deriter}
Let $\{\zeta_k\}$ and $\delta=\delta_0$ be defined as in Lemma
\ref{ineq}, with $w\in\mathcal{E}_{n}^{|J|}(\delta_0)$ such that
$\|h^o\|_{\infty}\leq 1$. If $|w'(\zeta)|\leq|\zeta||\log
|\zeta||^{|J|}$ in $D_{r+\frac{n-1}{n},\delta_0}$ then, there is a
positive constant $\delta_1$ such that, for  every $k\in\N$ and for
every $\zeta\in\overline{D}_{r+\frac{n-1}{n},\delta_1}$, we have the
inequality
\begin{equation}\label{ndz_h/dz}
\bigg|\frac{d\zeta_k}{d\zeta}\bigg|\leq
c\bigg|\frac{\zeta_k}{\zeta}\bigg|^{r+1+\frac{n-1}{n}}
\bigg|\frac{\log \zeta_k}{\log \zeta}\bigg|^{\frac{n-1}{n}},
\end{equation}
for some positive constant $c$.
\end{lemma}

\begin{proof}
Let us consider only the case $n\geq 3$, because if $n=2$ the
computations are easier. By \eqref{xshiftn} we get
\begin{displaymath}
\begin{aligned}
\frac{1}{\zeta_1^{r+\frac{n-1}{n}}(\log \zeta_1)^{\frac{n-1}{n}}}&=
\frac{1}{\zeta^{r+\frac{n-1}{n}}(\log
\zeta)^{\frac{n-1}{n}}}+r+\frac{n-1}{n}+
\left (\frac{n-1}{n}\right )\frac{1}{\log \zeta}\\
&\quad +O\left (\zeta^{\frac{1}{n}}(\log \zeta)^{i-\frac{n-1}{n}},
\frac{w(\zeta)}{\zeta^{\frac{1}{n}}(\log \zeta)^{\frac{1}{n}}},
\frac{w(\zeta)^{p-1}}{\zeta^{\frac{p-1}{n}}(\log
\zeta)^{\frac{p-1}{n}}} \right ),
\end{aligned}
\end{displaymath}
with $p$ integer and $3\leq p\leq n$. Our aim is to look for a
function $g(\zeta)$ such that
\begin{equation}\label{gn}
\begin{aligned}
\frac{1}{\zeta_1^{r+\frac{n-1}{n}}(\log
\zeta_1)^{\frac{n-1}{n}}}+g(\zeta_1)&=\frac{1}{\zeta^{r+\frac{n-1}{n}}
(\log \zeta)^{\frac{n-1}{n}}}+r+\frac{n-1}{n}+g(\zeta)\\
&\quad + O\left (\zeta^{\frac{1}{n}}(\log \zeta)^{i-\frac{n-1}{n}},
\frac{w(\zeta)}{\zeta^{\frac{1}{n}}(\log \zeta)^{\frac{1}{n}}},
\frac{w(\zeta)^{p-1}}{\zeta^{\frac{p-1}{n}}(\log
\zeta)^{\frac{p-1}{n}}} \right ).
\end{aligned}
\end{equation}
Comparing the previous equations we get
\begin{displaymath}
g(\zeta_1)-g(\zeta)=\left (\frac{1-n}{n}\right )\frac{1}{\log
\zeta}+ O\left (\zeta^{\frac{1}{n}}(\log \zeta)^{i-\frac{n-1}{n}},
\frac{w(\zeta)}{\zeta^{\frac{1}{n}}(\log \zeta)^{\frac{1}{n}}},
\frac{w(\zeta)^{p-1}}{\zeta^{\frac{p-1}{n}}(\log
\zeta)^{\frac{p-1}{n}}} \right ).
\end{displaymath}
The hypotheses on $w\in \mathcal{E}_{n}^{|J|}(\delta)$ and the
approximation
$g(\zeta_1)=g(\zeta)+g'(\zeta)(\zeta_1-\zeta)\\+O((\zeta_1-\zeta)^2)$
imply that
\begin{displaymath}
g'(\zeta)(\zeta_1-\zeta)+O((\zeta_1-\zeta)^2)= -\left
(\frac{n-1}{n}\right )\frac{1}{\log
\zeta}+O\left(\zeta^{\frac{1}{n}}(\log
\zeta)^{i-\frac{n-1}{n}}\right).
\end{displaymath}
From \eqref{xshiftn} we deduce that
\begin{displaymath}
\zeta_1-\zeta=-\zeta^{r+1+\frac{n-1}{n}}(\log
\zeta)^{\frac{n-1}{n}}+\zeta^{r+1+\frac{n-2}{n}}\phi(\zeta,w),
\end{displaymath}
with
$\zeta^{r+1+\frac{n-2}{n}}\phi(\zeta,w(\zeta))=O(\zeta^{r+2}(\log
\zeta)^i)$ and so we obtain
\begin{displaymath}
\begin{aligned}
-\zeta^{r+1+\frac{n-1}{n}}(\log \zeta)^{\frac{n-1}{n}}g'(\zeta)&+
O(\zeta^{r+2}(\log \zeta)^i)g'(\zeta)\\
&=-\left (\frac{n-1}{n}\right )\frac{1}{\log \zeta}+O
\left(\zeta^{\frac{1}{n}}(\log \zeta)^{i-\frac{n-1}{n}}\right),
\end{aligned}
\end{displaymath}
and hence
\begin{equation}\label{g'n}
g'(\zeta)=\left ( \frac{n-1}{n}\right )\frac{1}{\zeta^{r+1+
\frac{n-1}{n}}(\log \zeta)^{1+\frac{n-1}{n}}}.
\end{equation}

Let $g(\zeta)$ be a solution of the previous differential equation.
Then \eqref{gn} can be written as
\begin{equation}\label{varphin}
\frac{1}{\zeta_1^{r+\frac{n-1}{n}}(\log
\zeta_1)^{\frac{n-1}{n}}}+g(\zeta_1)=\frac{1}{\zeta^{r+\frac{n-1}{n}}
(\log \zeta)^{\frac{n-1}{n}}}+r+\frac{n-1}{n}+g(\zeta)+\varphi
(\zeta,w(\zeta)),
\end{equation}
where $\varphi$ is a function holomorphic in $\{\zeta,w,\log
\zeta\}$ and
\begin{displaymath}
\varphi (\zeta,w(\zeta))= O\left (\zeta^{\frac{1}{n}}(\log
\zeta)^{i-\frac{n-1}{n}}, \frac{w(\zeta)}{\zeta^{\frac{1}{n}}(\log
\zeta)^{\frac{1}{n}}}, \frac{w(\zeta)^{p-1}}{\zeta^{\frac{p-1}{n}}
(\log \zeta)^{\frac{p-1}{n}}} \right ).
\end{displaymath}

Now we want to estimate $\frac{d}{d\zeta}\varphi(\zeta,w(\zeta))$.
The hypotheses on $w(\zeta)$ and $w'(\zeta)$ imply that
\begin{displaymath}
\begin{aligned}
\frac{d}{d\zeta}\left (\frac{w(\zeta)}{\zeta^{\frac{1}{n}}(\log
\zeta)^{\frac{1}{n}}} \right )&=O\left
(\frac{w'(\zeta)}{\zeta^{\frac{1}{n}}(\log \zeta)^{\frac{1}{n}}},
\frac{w(\zeta)}{\zeta^{\frac{1}{n}+1}(\log
\zeta)^{\frac{1}{n}}}\right )=O\left(\zeta^{\frac{n-1}{n}}(\log
\zeta)^{|J|
-\frac{1}{n}}\right),\\
\frac{d}{d\zeta}\left (\zeta^{\frac{1}{n}}(\log
\zeta)^{i-\frac{n-1}{n}}\right )&=O\left
(\zeta^{\frac{1}{n}-1}(\log \zeta)^{i-\frac{n-1}{n}} \right ),\\
\frac{d}{d\zeta}\left
(\frac{w(\zeta)^{p-1}}{\zeta^{\frac{p-1}{n}}(\log
\zeta)^{\frac{p-1}{n}}} \right )&=O\left
(\frac{w'(\zeta)w(\zeta)^{p-2}}{\zeta^{\frac{p-1}{n}}(\log
\zeta)^{\frac{p-1}{n}}},\frac{w(\zeta)^{p-1}}{\zeta^{\frac{p-1}{n}+1}
(\log \zeta)^{\frac{p-1}{n}}}\right)\\
&=O\left( \zeta^{2(p-1)-1-\frac{p-1}{n}}(\log
\zeta)^{|J|(p-1)-\frac{p-1}{n}} \right).
\end{aligned}
\end{displaymath}
Therefore
\begin{equation}\label{varphi'n}
\frac{d}{d\zeta}\varphi(\zeta,w(\zeta))= O\left
(\zeta^{\frac{1}{n}-1}(\log \zeta)^{i-\frac{n-1}{n}} \right ).
\end{equation}

Iterating \eqref{varphin} we have
\begin{displaymath}
\begin{aligned}
\frac{1}{\zeta_k^{r+\frac{n-1}{n}}(\log
\zeta_k)^{\frac{n-1}{n}}}+g(\zeta_k)=\frac{1}{\zeta^{r+\frac{n-1}{n}}(\log
\zeta)^{\frac{n-1}{n}}}&+k\left (r+\frac{n-1}{n} \right
)+g(\zeta)\\
&+\sum_{j=0}^{k-1}\varphi (\zeta_j,w(\zeta_j)),
\end{aligned}
\end{displaymath}
and differentiating in $\zeta$ this last relation, we get
\begin{displaymath}
\begin{aligned}
-\bigg \{ \frac{(r+\frac{n-1}{n})\log \zeta_k +
\frac{n-1}{n}}{\zeta_k^{r+1+\frac{n-1}{n}}(\log
\zeta_k)^{1+\frac{n-1}{n}}} -g'(\zeta_k) \bigg
\}\frac{d\zeta_k}{d\zeta}&=-\bigg \{\frac{(r+\frac{n-1}{n})\log
\zeta +\frac{n-1}{n}}{\zeta^{r+1+\frac{n-1}{n}}
(\log \zeta)^{1+\frac{n-1}{n}}} -g'(\zeta)\bigg \}\\
&\quad
+\sum_{j=0}^{k-1}\frac{d}{d\zeta_j}\varphi(\zeta_j,w(\zeta_j))
\frac{d\zeta_j}{d\zeta}.
\end{aligned}
\end{displaymath}
Replacing \eqref{g'n} in the previous expression, we have
\begin{equation}\label{dz_h/dz}
\begin{aligned}
-\left (r+\frac{n-1}{n}\right )\frac{1}{\zeta_k^{r+1+\frac{n-1}{n}}
(\log \zeta_k)^{\frac{n-1}{n}}}\frac{d\zeta_k}{d\zeta}&=-\left
(r+\frac{n-1}{n}\right)\frac{1}{\zeta^{r+1+\frac{n-1}{n}}
(\log \zeta)^{\frac{n-1}{n}}}\\
&\quad
+\sum_{j=0}^{k-1}\frac{d}{d\zeta_j}\varphi(\zeta_j,w(\zeta_j))
\frac{d\zeta_j}{d\zeta}.
\end{aligned}
\end{equation}

At this point we can proceed by induction on $k$. If $k=1$, then
\eqref{dz_h/dz} becomes:
\begin{displaymath}
\frac{d\zeta_1}{d\zeta}=\frac{\zeta_1^{r+1+\frac{n-1}{n}}(\log
\zeta_1)^{\frac{n-1}{n}}}{\zeta^{r+1+\frac{n-1}{n}}(\log
\zeta)^{\frac{n-1}{n}}}\left \{
1-\frac{n}{nr+n-1}\zeta^{r+1+\frac{n-1}{n}}(\log
\zeta)^{\frac{n-1}{n}}\frac{d}{d\zeta}\varphi(\zeta,w(\zeta))
\right\},
\end{displaymath}
and from \eqref{varphi'n} we deduce that
\begin{displaymath}
\zeta^{r+1+\frac{n-1}{n}}(\log
\zeta)^{\frac{n-1}{n}}\frac{d}{d\zeta}\varphi(\zeta,w(\zeta))=
O(\zeta^{r+1}(\log \zeta)^i).
\end{displaymath}
For $\delta$ small enough, the quantities inside the brackets are
close to $1$, so there exists a constant $c$ such that
\begin{displaymath}
\bigg|\frac{d\zeta_1}{d\zeta}\bigg|\leq
c\bigg|\frac{\zeta_1}{\zeta}\bigg|^{r+1+\frac{n-1}{n}}\bigg|
\frac{\log \zeta_k}{\log \zeta}\bigg|^{\frac{n-1}{n}}.
\end{displaymath}
Let $k>1$. We can write \eqref{dz_h/dz} as
\begin{equation}\label{induz n}
\begin{aligned}
\frac{d\zeta_k}{d\zeta}&=\left \{ 1-\left
(r+\frac{n-1}{n}\right)^{-1}\zeta^{r+1+\frac{n-1}{n}} (\log
\zeta)^{\frac{n-1}{n}}\left(\sum_{j=0}^{k-1}\frac{d}{d\zeta_j}
\varphi(\zeta_j,w(\zeta_j))\frac{d\zeta_j}{d\zeta}\right)\right\}\\
&\qquad \cdot\frac{\zeta_k^{r+1+\frac{n-1}{n}}(\log
\zeta_k)^{\frac{n-1}{n}}}{\zeta^{r+1+\frac{n-1}{n}} (\log
\zeta)^{\frac{n-1}{n}}}.
\end{aligned}
\end{equation}
Using the inductive hypothesis, \eqref{varphi'n}, and Corollary
\ref{corserie} we get
\begin{displaymath}
\begin{aligned}
\sum_{j=0}^{k-1} \bigg |
\frac{d}{d\zeta_j}\varphi(\zeta_j,w(\zeta_j)) \bigg | \bigg |
\frac{d\zeta_j}{d\zeta} \bigg |&\leq
\frac{C}{|\zeta|^{r+1+\frac{n-1}{n}}|\log
\zeta|^{\frac{n-1}{n}}}\sum_{j=0}^{k-1}|\zeta_j|^{r+1}
|\log \zeta_j|^i\\
&\leq \frac{C}{|\zeta|^{r+1+\frac{n-1}{n}}|\log
\zeta|^{\frac{n-1}{n}}}\sum_{j=0}^{\infty}|\zeta_j|^{r+1}
|\log \zeta_j|^i \\
&\leq C_1 \frac{|\zeta|^{\frac{1}{n}}|\log
|\zeta||^{i-\frac{n-1}{n}}}{|\zeta|^{r+1+\frac{n-1}{n}}|\log
\zeta|^{\frac{n-1}{n}}},
\end{aligned}
\end{displaymath}
for some constant $C$ and $C_1$. Hence
\begin{displaymath}
\zeta^{r+1+\frac{n-1}{n}}(\log \zeta)^{\frac{n-1}{n}}\left
(\sum_{j=0}^{k-1}\frac{d}{d\zeta_j}
\varphi(\zeta_j,w(\zeta_j))\frac{d\zeta_j}{d\zeta}\right)=
O\left(\zeta^{\frac{1}{n}}(\log \zeta)^{i-\frac{n-1}{n}}\right)
\end{displaymath}
and the quantities inside the brackets in \eqref{induz n} are
uniformly close to $1$, so we get \eqref{ndz_h/dz}.
\end{proof}

\begin{lemma}
Let $T$ be defined as in \eqref{Tn}, and let $w\in
\mathcal{E}_{n}^{|J|}(\delta)$ with $\|h^o\|_{\infty}\leq 1$. If
$|w'(\zeta)|\leq|\zeta||\log |\zeta||^{|J|}$ for all $\zeta\in
D_{r+\frac{n-1}{n},\delta}$, then for $\delta$ small enough, we have
\begin{displaymath}
\bigg | \frac{d}{d\zeta}(Tw)(\zeta)\bigg |\leq|\zeta| |\log
|\zeta||^{|J|}, \quad \hbox{for all} \,\,\zeta\in
D_{r+\frac{n-1}{n},\delta}.
\end{displaymath}
\end{lemma}

\begin{proof}
Let us write
\begin{displaymath}
\frac{d}{d\zeta}(Tw)(\zeta)=\mathcal{S}_1+\mathcal{S}_2+\mathcal{S}_3.
\end{displaymath}
We have to estimate the three terms $\mathcal{S}_1$,
$\mathcal{S}_2$, and $\mathcal{S}_3$.

The term $\mathcal{S}_1$ comes from the differentiation of the term
$\zeta^{\frac{1}{n}}$, i.e.
\begin{displaymath}
\mathcal{S}_1=\frac{\zeta^{\frac{1}{n}-1}}{n} \sum_{k=0}^\infty
\zeta_k^{-\frac{1}{n}}H(\zeta_k,w(\zeta_k)).
\end{displaymath}
Using  \eqref{OHn} and Corollary \ref{corserie}, we get
\begin{displaymath}
|\mathcal{S}_1|\leq
\frac{1}{n|\zeta|^{\frac{n-1}{n}}}\sum_{k=0}^\infty|
\zeta_k|^{r+2+\frac{n-2}{n}}|\log \zeta_k|^{|J|}\leq C_1
|\zeta||\log |\zeta||^{|J|-\frac{n-1}{n}}\leq
\frac{1}{3}|\zeta||\log |\zeta||^{|J|}.
\end{displaymath}

The term $\mathcal{S}_2$ is given by
\begin{displaymath}
\mathcal{S}_2=\zeta^{\frac{1}{n}}\sum_{k=0}^\infty
\zeta_k^{-\frac{1}{n}} \frac{\partial}{\partial
w}H(\zeta_k,w(\zeta_k))\frac{dw}{d\zeta_k}\frac{d\zeta_k}{d\zeta}.
\end{displaymath}
From \eqref{OHn} and from the hypothesis on $w$, we have the
following estimate
\begin{displaymath}
\bigg |\frac{\partial}{\partial w}H(\zeta,w(\zeta)) \bigg | \leq
C_2\frac{|\zeta|^{r+\frac{n-1}{n}}}{|\log \zeta|^{\frac{1}{n}}},
\end{displaymath}
for some constant $C_2$. Using this inequality, Lemma \ref{deriter},
and the hypothesis on $w'(\zeta)$, we get
\begin{displaymath}
\begin{aligned}
\bigg |\frac{\partial}{\partial w}H(\zeta_k,w(\zeta_k))
\frac{dw}{d\zeta_k}\frac{d\zeta_k}{d\zeta}\bigg |& \leq C_3
\frac{|\zeta_k|^{r+\frac{n-1}{n}}}{|\log \zeta_k|^{\frac{1}{n}}}
|\zeta_k||\log |\zeta_k||^{|J|}\bigg|\frac{\zeta_k}
{\zeta}\bigg|^{r+1+\frac{n-1}{n}}\bigg|\frac{\log \zeta_k}{\log
\zeta}\bigg|^{\frac{n-1}{n}}\\
&\leq C_4\frac{|\zeta_k|^{2r+2+2\left(\frac{n-1}{n}\right)}|\log
|\zeta_k||^{|J| +\frac{n-2}{n}}}{|\zeta|^{r+1+\frac{n-1}{n}}|\log
|\zeta||^{\frac{n-1}{n}}},
\end{aligned}
\end{displaymath}
for some constants $C_3$ and $C_4$. Hence, using Corollary
\ref{corserie}, we obtain
\begin{displaymath}
\begin{aligned}
|\mathcal{S}_2|&\leq \frac{C_4}{|\zeta|^{r+1+\frac{n-2}{n}}|\log
|\zeta||^{\frac{n-1}{n}}}\sum_{k=0}^\infty
|\zeta_k|^{2r+2+2\left(\frac{n-1}{n}
\right)-\frac{1}{n}}|\log |\zeta_k||^{|J|+\frac{n-2}{n}}\\
&\leq C_5|\zeta||\log |\zeta||^{|J| -1}\leq \frac{1}{3}|\zeta||\log
|\zeta||^{|J|},
\end{aligned}
\end{displaymath}
for some constant $C_5$.

The term $\mathcal{S}_3$ is given by
\begin{displaymath}
\mathcal{S}_3=\zeta^{\frac{1}{n}}\sum_{k=0}^{\infty}
\frac{\partial}{\partial \zeta_k}\left
[\zeta_k^{-\frac{1}{n}}H(\zeta_k,w(\zeta_k)) \right
]\frac{d\zeta_k}{d\zeta}.
\end{displaymath}
Let us define $u(\cdot)$ by the equality
\begin{displaymath}
\mathcal{S}_3=\zeta^{\frac{1}{n}}\sum_{k=0}^{\infty}
\zeta_k^{-\frac{1}{n}}u(\zeta_k).
\end{displaymath}
An easy computation shows that
\begin{displaymath}
u(\zeta_k)=\left [-\frac{\zeta_k^{-1}}{n}H(\zeta_k,w(\zeta_k))+
\frac {\partial}{\partial \zeta_k}H(\zeta_k,w(\zeta_k)) \right
]\frac{d\zeta_k}{d\zeta}.
\end{displaymath}
Since $H(\zeta,w(\zeta))=O\left(\zeta^{r+2+\frac{n-1}{n}}(\log
\zeta)^{|J|}\right)$, we get
\begin{displaymath}
\frac{\partial}{\partial
\zeta}H(\zeta,w(\zeta))=O\left(\zeta^{r+1+\frac{n-1}{n}} (\log
\zeta)^{|J|}\right),
\end{displaymath}
and hence
\begin{displaymath}
|u(\zeta_k)|\leq C_6 |\zeta_k|^{r+1+\frac{n-1}{n}} |\log
\zeta_k|^{|J|}\bigg |\frac{d\zeta_k}{d\zeta}\bigg |,
\end{displaymath}
for some constant $C_6$ independent of $w$ and $w'$. Using Lemma
\ref{deriter}, we can proceed to a further estimate
\begin{displaymath}
\begin{aligned}
|u(\zeta_k)|&\leq \frac{C_7}{|\zeta|^{r+1+\frac{n-1}{n}}|\log
\zeta|^{\frac{n-1}{n}}}|\zeta_k|^{2r+2+2\left(
\frac{n-1}{n}\right)}|\log \zeta_k|^{|J|+\frac{n-1}{n}}.
\end{aligned}
\end{displaymath}
Finally, using  Corollary \ref{corserie} again, there exists a
constant $C_8$ such that
\begin{displaymath}
\begin{aligned}
|\mathcal{S}_3|&\leq \frac{C_7}{|\zeta|^{r+1+\frac{n-2}{n}}|\log
\zeta|^{\frac{n-1}{n}}}\sum_{k=0}^{\infty}|\zeta_k|^{2r+2+2\left(
\frac{n-1}{n}\right)-\frac{1}{n}}|\log \zeta_k|^{|J|+\frac{n-1}{n}}\\
&\leq C_8|\zeta||\log |\zeta||^{|J|-\frac{n-1}{n}}\leq \frac{1}{3}
|\zeta||\log |\zeta||^{|J|}.
\end{aligned}
\end{displaymath}
Therefore, for $\delta$ small enough, we get, for all $\zeta\in
D_{r+\frac{n-1}{n},\delta}$,
\begin{displaymath}
\bigg | \frac{d}{d\zeta}(Tw)(\zeta)\bigg |\leq
|\mathcal{S}_1|+|\mathcal{S}_2|+|\mathcal{S}_3|\leq|\zeta||\log
|\zeta||^{|J|}.
\end{displaymath}
\end{proof}

From this last Lemma it follows that the operator $T$ sends the
convex closed set $\mathcal{F}_{n}^{|J|}(\delta)$ into itself. It
remains to prove that $T$ is a contraction. We need a last result.

\begin{lemma}\label{diffzn}
Let $\hat{f}$ be a map of the form \eqref{xshift} or
\eqref{xshiftn}, with $h=2n-3$. Let $u(\zeta)=\zeta^{2}(\log
\zeta)^{|J|} h_1(\zeta)$ and $v(\zeta)=\zeta^{2}(\log \zeta)^{|J|}
h_2(\zeta)$ be two functions in $\mathcal{F}_{n}^{|J|}(\delta)$. Let
$\{\zeta_k\}$ and $\{\zeta'_k\}$ be, respectively, the iterates of
$\zeta$ by $\hat{f}_1(\zeta,u(\zeta))$ and
$\hat{f}_1(\zeta,v(\zeta))$. Then, for $\delta$ small enough, there
exists a constant $K$ such that, for all $k\in\mathbb N$
\begin{displaymath}
|\zeta'_k-\zeta_k|\leq  K|\zeta|^{3-\frac{1}{n}}|\log |\zeta||^{|J|
-\frac{1}{n}}\|h_2-h_1\|_\infty.
\end{displaymath}
\end{lemma}

\begin{proof}
For $\zeta$ and $\zeta'$ in the same connected component of
$D_{r+\frac{n-1}{n},\delta}$, we estimate the variation
\begin{displaymath}
\hat{f}_1(\zeta',v(\zeta'))-\hat{f}_1(\zeta,u(\zeta)).
\end{displaymath}
From \eqref{xshift} or \eqref{xshiftn}, with $h=2n-3$, we have
\begin{equation}\label{varn}
\begin{aligned}
&\hat{f}_1(\zeta',v(\zeta'))-\hat{f}_1(\zeta,u(\zeta))\\
&=\left [1-\frac{(\zeta')^{r+1+\frac{n-1}{n}}(\log
\zeta{'})^{\frac{n-1}{n}} -\zeta^{r+1+\frac{n-1}{n}}(\log
\zeta)^{\frac{n-1}{n}}}{\zeta'-\zeta}+O(|\zeta''|^{r+1}
|\log \zeta''|^i)\right ] \\
&\quad \cdot (\zeta'-\zeta)+(v(\zeta')-u(\zeta))
O\left(|\zeta''|^{r+1+\frac{n-2}{n}}|\log
\zeta''|^{\frac{n-2}{n}}\right),
\end{aligned}
\end{equation}
where $|\zeta''|=\max\{|\zeta|,|\zeta'|\}$. We observe that when $k$
goes to infinity then $\frac{\zeta'_k}{\zeta_k} \sim 1$ by
\eqref{inequalities}, so we can replace $|\zeta''|$ by $|\zeta|$ in
these estimates. Therefore
\begin{equation}\label{diffn}
\begin{aligned}
v(\zeta')-u(\zeta)&=v(\zeta')-v(\zeta)+v(\zeta)-u(\zeta)\\
&=O(|\zeta||\log |\zeta||^{|J|})(\zeta'-\zeta)+O(|\zeta|^2|\log
|\zeta||^{|J|})\|h_2-h_1\|_\infty.
\end{aligned}
\end{equation}
Using this last relation, \eqref{varn} becomes
\begin{displaymath}
\begin{aligned}
\zeta'_1-\zeta_1=&(\zeta'-\zeta)\bigg
[1-\frac{(\zeta')^{r+1+\frac{n-1}{n}}(\log \zeta{'})^{\frac{n-1}{n}}
-\zeta^{r+1+\frac{n-1}{n}}(\log \zeta)^{\frac{n-1}{n}}}{\zeta'-\zeta} \\
&\qquad \qquad + O(|\zeta|^{r+1}|\log |\zeta||^i)\bigg ] \\
&+O\left(|\zeta|^{r+3+\frac{n-2}{n}}|\log |\zeta||^{|J|
+\frac{n-2}{n}}\right)\|h_2-h_1\|_\infty.
\end{aligned}
\end{displaymath}
For $\zeta$ and $\zeta'\in D_{r+\frac{n-1}{n},\delta}$, and $\delta$
small enough, the modulus of the quantity inside the brackets is
bounded by $1$, and so we get
\begin{displaymath}
|\zeta'_1-\zeta_1|\leq |\zeta'-\zeta| +
K_1|\zeta|^{r+3+\frac{n-2}{n}} |\log
|\zeta||^{|J|+\frac{n-2}{n}}\|h_2-h_1\|_\infty,
\end{displaymath}
for some constant $K_1$. Iterating this process and taking
$\zeta=\zeta'$, we get, for all $k\in \N$
\begin{displaymath}
|\zeta'_k-\zeta_k|\leq K_1\sum_{p=0}^\infty
|\zeta_p|^{r+3+\frac{n-2}{n}}|\log |\zeta_p||^{|J|
+\frac{n-2}{n}}\|h_2-h_1\|_\infty
\end{displaymath}
and by Corollary \ref{corserie}
\begin{displaymath}
|\zeta'_k-\zeta_k|\leq K|\zeta|^{3-\frac{1}{n}}|\log |\zeta||^{|J|
-\frac{1}{n}}\|h_2-h_1\|_\infty,
\end{displaymath}
for some constant $K$.
\end{proof}

\begin{prop}
Let $T$ be defined as in \eqref{Tn}. For $\delta$ small enough, the
restriction of $T$ is a contraction of the convex closed subset
$\mathcal{F}_{n}^{|J|}(\delta)$ of the Banach space
$\mathcal{E}_{n}^{|J|}(\delta)$ into itself.
\end{prop}

\begin{proof}
Let $u(\zeta)=\zeta^2(\log \zeta)^{|J|} h_1(\zeta)$ and
$v(\zeta)=\zeta^2(\log \zeta)^{|J|} h_2(\zeta)$ be two functions in
$\mathcal{F}_{n}^{|J|}(\delta)$. We can write the variation
\begin{displaymath}
Tu(\zeta)-Tv(\zeta)=\zeta^{\frac{1}{n}}\sum_{k=0}^\infty
\left[\zeta_k^{-\frac{1}{n}}H(\zeta_k,u(\zeta_k))-
(\zeta'_k)^{-\frac{1}{n}}H(\zeta'_k,v(\zeta'_k))\right]
\end{displaymath}
as the sum of two terms
\begin{displaymath}
\mathcal{R}_1=\zeta^{\frac{1}{n}}\sum_{k=0}^\infty
\zeta_k^{-\frac{1}{n}}[H(\zeta_k,u(\zeta_k))-H(\zeta'_k,v(\zeta'_k))]
\end{displaymath}
and
\begin{displaymath}
\mathcal{R}_2=\zeta^{\frac{1}{n}}\sum_{k=0}^\infty
\left[\zeta_k^{-\frac{1}{n}}-(\zeta'_k)^{-\frac{1}{n}}\right]
H(\zeta'_k,v(\zeta'_k)).
\end{displaymath}
From \eqref{OHn} and \eqref{inequalities} there exists a constant
$K_1$ such that
\begin{displaymath}
\begin{aligned}
&|H(\zeta_k,u(\zeta_k))-H(\zeta'_k,v(\zeta'_k))|\\
&\leq K_1 \left[\frac{|\zeta_k|^{r+\frac{n-1}{n}}}{|\log
|\zeta_k||^{\frac{1}{n}}}|u(\zeta_k)-v(\zeta'_k)| +
|\zeta_k|^{r+1+\frac{n-1}{n}}|\log |\zeta_k||^{|J|}
|\zeta_k-\zeta'_k|\right].
\end{aligned}
\end{displaymath}
On the other hand, using \eqref{diffn} we can write
\begin{displaymath}
|v(\zeta'_k)-u(\zeta_k)|\leq |\zeta_k||\log |\zeta_k||^{|J|}
|\zeta'_k-\zeta_k| + |\zeta_k|^2|\log
|\zeta_k||^{|J|}\|h_2-h_1\|_\infty,
\end{displaymath}
and by Lemma \ref{diffzn} it follows that
\begin{displaymath}
\begin{aligned}
|v(\zeta'_k)-u(\zeta_k)|\leq &\left(K|\zeta_k||\log
|\zeta_k||^{|J|}|\zeta|^{3-\frac{1}{n}}|\log |\zeta||^{|J|
-\frac{1}{n}} + |\zeta_k|^2|\log |\zeta_k||^{|J|}
\right)\\
&\cdot \|h_2-h_1\|_\infty.
\end{aligned}
\end{displaymath}
Applying Corollary \ref{corserie}, we get
\begin{displaymath}
|\mathcal{R}_1|\leq K_7|\zeta|^2|\log |\zeta||^{|J| -1}
\|h_2-h_1\|_\infty,
\end{displaymath}
for some constant $K_7$. In fact, there exist constants $K_i$, with
$i=2,\ldots,7$ such that
\begin{displaymath}
\begin{aligned}
|\mathcal{R}_1|&\leq K_1|\zeta|^{\frac{1}{n}}\\
&\cdot \sum_{k=0}^\infty \frac{1}{|\zeta_k|^{\frac{1}{n}}}\left|
\frac{|\zeta_k|^{r+\frac{n-1}{n}}}{|\log
|\zeta_k||^{\frac{1}{n}}}|u(\zeta_k)-v(\zeta'_k)| +
|\zeta_k|^{r+1+\frac{n-1}{n}}|\log |\zeta_k||^{|J|}
|\zeta_k-\zeta'_k|\right|\\
&\leq K_2|\zeta|^{3}|\log |\zeta||^{|J|
-\frac{1}{n}}\|h_2-h_1\|_\infty\sum_{k=0}^\infty
|\zeta_k|^{r+1+\frac{n-2}{n}} |\log
|\zeta_k||^{|J| -\frac{1}{n}}\\
&\quad + K_1|\zeta|^{\frac{1}{n}}\|h_2-h_1\|_\infty\sum_{k=0}^\infty
|\zeta_k|^{r+2+\frac{n-2}{n}} |\log |\zeta_k||^{|J| -\frac{1}{n}}\\
&\quad + K_3|\zeta|^{3}|\log
|\zeta||^{|J|-\frac{1}{n}}\|h_2-h_1\|_\infty\sum_{k=0}^\infty
|\zeta_k|^{r+1+\frac{n-2}{n}}|\log
|\zeta_k||^{|J|}\\
&\leq K_4|\zeta|^{4-\frac{1}{n}}|\log |\zeta||^{2|J|
-1-\frac{1}{n}}\|h_2-h_1\|_\infty
+ K_5|\zeta|^2|\log |\zeta||^{|J| -1}\|h_2-h_1\|_\infty\\
&\quad + K_6|\zeta|^{4-\frac{1}{n}}|\log
|\zeta||^{2|J|-1}\|h_2-h_1\|_\infty\leq K_7 |\zeta|^2|\log
|\zeta||^{|J| -1}\|h_2-h_1\|_\infty.
\end{aligned}
\end{displaymath}

We can estimate the term $\mathcal{R}_2$ as follows: first we can
write
\begin{displaymath}
\left|\zeta_k^{-\frac{1}{n}}-(\zeta'_k)^{-\frac{1}{n}}\right|=
\left|\zeta_k^{-\frac{1}{n}}\right|\left|1-\exp
\left(-\frac{1}{n}\log \frac{\zeta'_k}{\zeta_k}\right)\right|,
\end{displaymath}
then by \eqref{OHn} and since $v\in \mathcal{F}_{n}^{|J|}(\delta)$,
it follows that
\begin{displaymath}
\left|1-\exp \left(-\frac{1}{n}\log
\frac{\zeta'_k}{\zeta_k}\right)\right||H(\zeta'_k,v(\zeta'_k))| \leq
K'_1 \frac{|\zeta'_k-\zeta_k|}{|\zeta_k|}
|\zeta_k|^{r+2+\frac{n-1}{n}}|\log |\zeta_k||^{|J|},
\end{displaymath}
and using Lemma \ref{diffzn}, we get
\begin{displaymath}
\begin{aligned}
&\left|1-\exp \left(-\frac{1}{n}\log
\frac{\zeta'_k}{\zeta_k}\right)\right||H(\zeta'_k,v(\zeta'_k))|\\
&\qquad \leq K'_2|\zeta|^{3-\frac{1}{n}}|\log |\zeta||^{|J|
-\frac{1}{n}}\|h_2-h_1\|_\infty|\zeta_k|^{r+1+\frac{n-1}{n}}|\log
|\zeta_k||^{|J|},
\end{aligned}
\end{displaymath}
where $K_1'$ and $K_2'$ are some constants. Finally by Corollary
\ref{corserie}
\begin{displaymath}
|\mathcal{R}_2|\leq K'_3|\zeta|^{4-\frac{1}{n}}|\log |\zeta||^{2|J|
-1}\|h_2-h_1\|_\infty,
\end{displaymath}
for some constant $K_3'$.

We observe that the term $\mathcal{R}_2$ is of smaller order than
$\mathcal{R}_1$. Hence, for $\delta$ small enough, there exists a
constant $K'$ such that
\begin{displaymath}
|Tu(\zeta)-Tv(\zeta)|\leq K'|\zeta|^{2}|\log |\zeta||^{|J|
-1}\|h_2-h_1\|_\infty\leq C|\zeta|^{2}|\log |\zeta||^{|J| }\
|h_2-h_1\|_\infty,
\end{displaymath}
for a suitable $C<1$. By definition of the norm in
$\mathcal{F}_{n}^{|J|}(\delta)$, this implies that, for $\delta$
small enough, we have
\begin{displaymath}
\|Tu-Tv\|\leq C\|u-v\|.
\end{displaymath}
So $T$ is a contraction.
\end{proof}

Since $n\geq 2$, we observe that $D_{r+\frac{n-1}{n},\delta}$ has
$r+1$ connected (and simply connected) components having $O$ on
their boundary. Hence, if $\hat{f}$ is defined by \eqref{fhat} or
\eqref{fhatn}, we obtain $r+1$ parabolic curves for $\hat{f}$ at the
origin. These parabolic curves are distinct as sets. We think that
the union of any two curves is not contained in a larger parabolic
curve, but for now we do not have a proof of this statement.
However, there exists at least one parabolic curve for $\hat{f}$ at
the origin and hence there exists at least one parabolic curve at
the origin for $f$ given by \eqref{fstartn}.

This completes the proof of Theorem \ref{m=n-1}, and hence of our
main result, Theorem ~\ref{teomain}.

\bibliographystyle{amsplain}

\end{document}